\documentclass[a4paper, 11pt]{article}
\usepackage{amsmath} \usepackage{amsfonts} \usepackage{amssymb}\usepackage{latexsym}
\usepackage[english]{babel}
\usepackage{amscd}
\usepackage[dvips,final]{graphics}
\usepackage{locengli,math}
\usepackage{graphicx,pst-plot,psfig,pstricks,pst-node,
pst-text,pst-tree,pst-char,pst-coil,Addaletree}
\newcommand{\N}{\hat{\mathbb{N}}}
\newcommand{\B}{\mathcal{B}}

\newcommand{\M}{\mathcal{M}}
\newcommand{\Na}{\mathbb{N}}
\newcommand{\Nb}{\mathbb{N}_b}
\def\YY{Y\!\!\!\! Y}

\def\NN{N\!\!\!\! N}

\def\u{\underline }
\newcommand{\Ng}{\hat{\NN}}
\newcommand{\DM}{ \textsf{Dend}_{\B}(\M)}
\newcommand{\TrM}{ \textsf{TriDend}_{\B}(\M)}
\newcommand{\NCP}{ \textsf{NCP}}
\newcommand{\Pro}{ \textsf{Pr}}
\newcommand{\otm}{\otimes_{ _{\B}}}

\begin{document}
\begin{center}
\textbf{\LARGE{\textsf{Free dendriform dialgebras: reformulation and application in free probability. Part I}}}
\footnote{Supported by the European Commission HPRN $\sim$ CTN2002 $\sim$ 00279, RTN QP-Applications.
}
\vskip1cm
\parbox[t]{14cm}{\large{
Philippe {\sc Leroux}}\\
\vskip4mm
{\footnotesize
\baselineskip=5mm
Institut f\"ur Mathematik und Informatik,\\
Ernst-Moritz-Arndt-Universit\"at, Jahnstra$\beta$e 15a, 17487 Greifswald, Germany,\\ leroux@uni-greifswald.de}}
\end{center}
\vskip1cm
\baselineskip=5mm
\noindent
\begin{center}
\begin{large}
\textbf{ 30/06/04}
\end{large}
\end{center}
{\bf Abstract:} 
We propose a reformulation of some results known on the free dendriform dialgebra on one generator from a parenthesis setting. This turns out to be more tractable and simplify proofs. We develop also the arithmetree on planar rooted binary trees and point out 
a connection to free probability by identifying noncrossing partitions with binary trees and by introducing the concept of \textsf{NCP}-operad. 
\section{Introduction}
In this paper, $K$ is a null characteristic field, $\Na$ is the semiring of integers and $\Nb^n$ stands for the set $\{\vec{v}:=(v_1, \ldots, v_n) \in \Na^n; \  \forall \ 1\leq i \leq n, \ \ 0< v_i \leq i \}$, in bijection with
$S_n$, the symmetric group over $n$ elements. If $S$ is a finite set, then $\card(S)$ denotes its cardinal,  $KS$, the $K$-vector space spanned by $S$ and $\bra S \ket$, the free associative semigroup
generated by $S$.
Rooted planar binary trees will be called binary trees for short. If $k \in K$, $\vec{v} \in K^n$, then $(k+ \vec{v})$ is the vector $\vec{v}$
whose coordinates have been shifted by $k$.
In Section 2, we recall briefly what regular operads mean. In Section 3, we propose a reformulation of the free dendriform dialgebra over the generator $\treeA$ \textit{via} a parenthesis setting. This framework has the advantage to make proofs easier. We propose in the same time both a brief survey on trees and new results proved from the parenthesis setting.
In Section 4, we present a bijection between planar rooted binary trees and noncrossing partitions. Noncrossing partitions can be viewed from
a rooted planar binary trees by `projecting' SW-NE branches on a particular axis.
This allows the introduction of the concept of \textsf{NCP}-operads, whose axioms look like regular operads ones. We conclude by proposing
a connection between free probability and the free dendriform dialgebra setting.
\section{Brief recall on $K$-linear regular operads}
We follow \cite{Lodayren}. Given a $K$-algebra $A$ `of type $\mathcal{P}$', one considers the family of the $K$-vector spaces $\mathcal{P}(n)$ of $n$-ary operations. Therefore, we have a linear map
$ \Phi: \mathcal{P}(n) \otimes A^{\otimes n} \xrightarrow{} A, \ \ \Phi(f; (a_1, \ldots, a_n)) \mapsto f(a_1, \ldots, a_n).$
Operations can be composed in the following natural ways. For $f \in \mathcal{P}(m)$, $g \in \mathcal{P}(n)$; $\forall \ 1 \leq i \leq m,$ $f \circ_i g \in \mathcal{P}(m+n-1)$ is defined by:
$$ f \circ_i g (a_1, \ldots, a_{m+n-1}):= f(a_1, \ldots,a_{i-1}, g(a_i, \ldots,a_{i+n-1}), a_{i+n}, \ldots a_{m+n-1}). $$
These composition operations have to obey natural conditions \cite{Lodayren} which are parenthesing compatibilities.
If $h \in \mathcal{P}(l)$, $f \in \mathcal{P}(m)$ and $g \in \mathcal{P}(n)$, then
$ (h \circ_i f) \circ_{j+m-1} g = (h \circ_j g) \circ_i f; \ \ 1\leq i < j\leq l, \ \
 (h \circ_i f) \circ_{i+j-1 } g = h \circ_i (f \circ_j g); \ \ 1 \leq i \leq l; \ 1 \leq j \leq m.$
A $K$-\textit{linear regular operad} $\mathcal{P}$ is then a family of $K$-vector spaces $(\mathcal{P}(n))_{n>0}$ equipped with composition maps $\circ_i$ verifying the above relations. If all possible operations are generated by composition from $\mathcal{P}(2)$, then the operad is said to be {\it{binary}}. It is said to be {\it{quadratic}} if all the relations between operations are consequences of relations described exclusively with the help of monomials with two operations. 
In this case, the free $\mathcal{P}$-algebra is entirely induced by the free $\mathcal{P}$-algebra on one generator $\mathcal{P}(K):= \oplus_{n \geq 1} \ \mathcal{P}(n)$. The generating function of the regular operad $\mathcal{P}$ is given by:
$ f^{\mathcal{P}}(x):= \sum \ (-1)^n \textrm{dim} \ \mathcal{P}(n) x^n.$
Below, we will indicate the sequence $(\textrm{dim} \ \mathcal{P}(n))_{n \geq 1}$.
Let $V$ be a $K$-vector space. The {\it{free $\mathcal{P}$-algebra}} $\mathcal{P}(V)$ on $V$ is by definition,
a $\mathcal{P}$-algebra equipped with a linear map $i: \ V \mapsto \mathcal{P}(V)$ which satisfies the following universal property:
for any linear map $f: V \xrightarrow{} A$, where $A$ is a $\mathcal{P}$-algebra, there exists a unique  $\mathcal{P}$-algebra morphism $\bar{f}: \mathcal{P}(V) \xrightarrow{} A$ such that $\bar{f} \circ i=f$.
Since our $\mathcal{P}$-algebras are regular, the free $\mathcal{P}$-algebra over a $K-$vctor space $V$ is of the form:
$ \mathcal{P}(V) := \bigoplus_{n \geq 1} \ \mathcal{P}_n \otimes V^{\otimes n} .$
In particular, the free $\mathcal{P}$-algebra on one generator $x$ is $ \mathcal{P}(Kx) := \bigoplus_{n \geq 1} \mathcal{P}_n$.
In the sequel, our two $K$-linear operads $\mathcal{P}$ will have only two binary operations $\bullet_1$, $\bullet_2$ generating $\mathcal{P}_2$ and three constraints in $\mathcal{P}_3$. Therefore, on one generator $x$,
$\mathcal{P}_1:= Kx$, $\mathcal{P}_2:= K(x \bullet_1 x) \oplus K(x\bullet_2 x)$.
The space of three variables made out of two operations is of dimension $2 \times 2^2= 8$. As we have three relations or constraints, the space $\mathcal{P}_3:=K((x \bullet_i x) \bullet_j x) \oplus K(x \bullet_i (x \bullet_j x)),$ for  $i,j:=1,2,$ has the dimension equal to $8-3=5.$
The sequence associated with the dimensions of $(\mathcal{P}_n)_{n \in \mathbb{N}}$ starts with $1, \ 2, \ 5, \ldots$, which
is the beginning of the Catalan numbers sequence.
\section{Arithmetics on trees from operads}
Dendriform dialgebras have been introduced by J.-L. Loday \cite{Loday}
as dual, in the operadic sense, to associative dialgebras, themselves motivated by $K$-theory. The free dendriform dialgebra on one generator is then closely related to binary trees. Major developments
have been put forward by using the Hopf algebra structure on the regular representations of the permutation groups found by C. Reutenauer and C. Malvenuto \cite{ReuMal} and connections between permutations and binary trees. Since then, an arithmetic
on trees have been introduced by J.-L. Loday \cite{Lodayarithm}.
The aim of this section is to present another way to handle the free dendriform dialgebra on one generator. Instead of starting with coding binary trees \textit{via} permutations, we focus on the parenthesing meaning of binary trees. In addition to be simpler, we hope this viewpoint will be more tractable for future computer developments. Another aim of this section is to put the arithmetics found in \cite{Lodayarithm} at the heart of the free dendriform dialgebra on one generator, to recover already known results with different and easier proofs, to produce extra results and at the same time to give a survey
of binary trees viewed from an operadic point of view.

A tree is binary if any vertex is trivalent. The set of planar rooted binary
trees with $n$ vertices, so called also $n$-trees, and considered up to isotopies, will be denoted by $Y_n$ (\textit{i.e.}, $n+1$ leaves and one root). The integer $n$ is also called the degree of a tree of $Y_n$, and $\card(Y_n) =c_n$, the Catalan numbers. In low dimensions, these sets are:
$$Y_0:=\{ \treeO \}, \ Y_1:=\{ Y:=\treeA \}, \ Y_2:=\{ \treeAB, \treeBA \}, \ Y_3:=\{ \treeABC, \treeBAC, \treeACA, \treeCAB, \treeCBA \}.$$
Arithmetree has been introduced by J.-L. Loday \cite{Lodayarithm} and is the analogue of the usual semiring $\Na$
at the level of planar binary trees. To present such binary trees, a code, rooted in permutations, has been introduced in \cite{Loday}, see also \cite{LRbruhat, Lodayarithm} and  is related to the grafting operation. The grafting of  a $p$-tree $\tau_1$ with a $q$-tree $\tau_2$ gives a $p+q+1$-tree denoted by $\tau_1 \vee \tau_2$ obtained by identifying the root of $\tau_1$ (resp.  $\tau_2$) with the left (resp. the right) leaf of $Y$. Set $[0]:= \treeO$. Any binary trees can be encoded into a sequence of integers by computing the rule $[\tau_1,p+q+1,\tau_2]$, where $\tau_1$, (resp. $\tau_2$) stands for the sequence of integers associated with the $p$-tree $\tau_1$, (resp. $\tau_2$). For instance, $\treeA:=[0] \vee [0] =[1]$, $\treeACA:=[1] \vee [1]:= [131]$. The sequence of integers associated with trees above are (from left to right):
$$ [0], [1], [12], [21], [123], [213], [131], [312], [321]. $$
This labelling has many advantages but do not fit well with the Tamari order of $Y_n$. Indeed, $Y_n$ can be endowed with a poset structure, often called the Tamari lattice,  by declaring that $\pi < \tau$, (also denoted by  $\pi \rightarrow \tau$) if $\tau$ can be obtained from $\pi$ by moving edges from left to right. For instance:
$ \treeAB  \rightarrow \treeBA.$
\subsection{Binary trees versus vectors}
In this subsection, we propose another way to encode binary trees which is compatible with the Tamari order. For that,
we associate with a planar binary tree of $Y_n$ a unique vector of $\Na^{n}$ in the following way. To any binary tree $\tau$ corresponds a unique parenthesing, and therefore a unique monomial
in $\bra \ x_1, \ldots, x_n, ( ,)\ \ket$ and thus  a unique monomial
in $\bra \ x_1, \ldots, x_n, ( \ \ket$ obtained by forgeting all right parentheses. Proceeding this way, we obtain an injection:
$ Exp:   Y_n \xrightarrow{} \bra \ x_1, \ldots, x_{n+1}, ( \ \ket.$
In the sequel, to ease notation, the unique parenthesing
associated with the binary tree $\tau$ will be also represented by $Exp(\tau)$ as in the following example. 
\begin{center}
\includegraphics*[width=7cm]{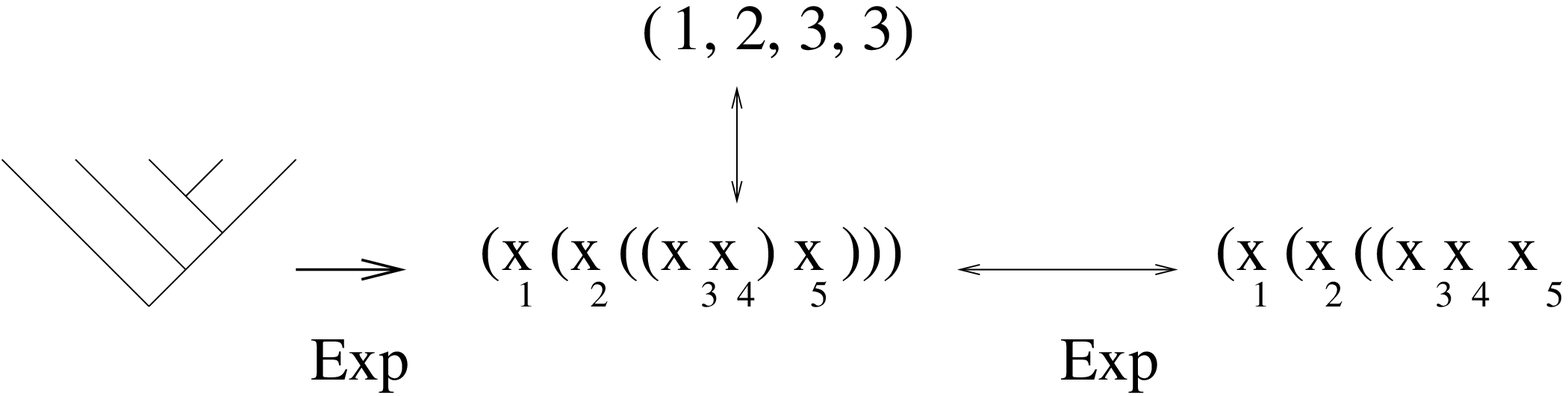}
\end{center}
Encode the parentheses of $Exp(\tau)$ of the binary tree $\tau$ in a vector  $\vec{v} :=(v_1, v_2, \ldots, v_{n})$ of $\Na^{n}$ by declaring that for all $1 \leq i \leq n$,
$v_i :=i$ if and only if there exists a left parenthesing at the left hand side of $x_i$, \textit{i.e.}, $\ldots(^p x_i \ldots$, with $p>0$, occurs in the monomial $Exp(\tau)$. Otherwise, there exists a unique most right parenthesis at the right hand side of $x_i$ which closes a unique left parenthesis say open at $x_j$. In this case, $v_i :=j$. Observe that this framework works since binary trees, \textit{via} their leaves, model all parentheses one can obtained from a  binary operation. We then obtain an injective map:
$ \vec{name}:   Y_n \xrightarrow{} \Na^{n},$
which map any tree $\tau$ into a vector, $\vec{name}(\tau)$, also denoted by $\vec{\tau}$ for short, called the name of $\tau$.
This coding can be extended to 
$ \Na^{n+1}$ by observing that coding the last leave corresponding to   $x_{n+1}$ gives always 1. In the sequel, $\vec{name}(Y_n)$ will be denoted by $ \N^{n}$ and by complete expression, we mean a monomial
of $\bra \ x_1, \ldots, x_n, ( ,)\ \ket$ in one-to-one correspondence with a rooted planar binary tree, \textit{i.e.,} every $($ is closed by a unique $)$.
\begin{prop}
\label{err}
Let  $\vec{v}  \in \Nb^n$. There exists a unique
monomial  
$(^{q_1}x_1, (^{q_2}x_2\ldots, (^{q_n}x_n x_{n+1}$ from $\bra \ x_1, \ldots, x_{n+1}, (, ) \ \ket$ associated with $\vec{v}$, where
$q_i$ is the number of $i$ appearing in $\vec{v}$. 
Such an algorithm gives a surjective map
$Tree: \Nb^n \xrightarrow{} Y_n$.
\end{prop}
\Proof
We proceed by induction.
Fix $\vec{v}:=(v_1, \ldots, v_n)  \in \Nb^n$. Its associated monomial in
$\bra \ x_1, \ldots, x_n, ( \ \ket$ is of the form $(^{q_1}x_1(^{q_2}x_2 \ldots (^{q_n}x_nx_{n+1}$, where $q_j$ is the number of $j$ appearing in $\vec{v}$. Observe that  $\sum q_j=n$ and $1 \leq q_j \leq n-j+1$ since $j$ may appear only from $x_j$. Take the highest $j$, with $q_j \not=0$, \textit{i.e.}, consider $(^{q_j}x_j x_{j+1} \ldots x_n$. As parentheses model a binary operation, there is a unique way to set right parentheses, namely $(^{q_j}x_jx_{j+1}), \ldots, x_{j + q_j}) \ldots x_nx_{n+1}$. This gives a complete expression $X:=(^{q_j}x_jx_{j+1}), \ldots, x_{j + q_j}) $ and  $(^{q_1}x_1(^{q_2}x_2 \ldots (^{q_{j-1}}x_{j-1}X \ldots x_{n+1}$ with $\sum q_k=n-q_j$. The proof is complete by induction.
\eproof

\noindent
Proposition \ref{err} is in fact a correcting error-code. Let us apply it to $(1,2,1,2)$. This gives $((x_1((x_2x_3x_4x_5$, \textit{i.e.}, $((x_1((x_2x_3)x_4))x_5)$, \textit{i.e.},
the tree named by $(1,2,2,1)$. 
\begin{coro}
(Reconstruction criterion)
A vector $\vec{v} \in \Nb^n$ is the name of a binary tree if and only if
$name(Tree(\vec{v}))=\vec{v}$.
\end{coro}
Another equivalent way to describe this coding consists to  start with $x_i$, go to the left and count the number of $x_k$, $k \leq i$, and the number of left parentheses. When these two numbers fit, take the last encountered $x_j$ and set  $v_i:=j$. This description can be found in B. E. Sagan \cite{Sagan}, whose one of his motivations was to compute the M\"obius function $M$ for the Tamari lattice. To state Theorem \ref{Sagan}, introduce the set $\mathcal{A}(L)$ of a lattice $L$ with minimal element, to denote the set of all atoms of $L$, ---those elements such that there is no other one between them and the minimum---. Such a set is call \textit{independent} if for all $B \subsetneq \mathcal{A}(L)$, $\bigvee B < \bigvee \mathcal{A}(L)$, where $\bigvee$ stands for the least upper bound operation. The following result holds. 
\begin{theo}[B. E. Sagan \cite{Sagan}]
\label{Sagan}
Let $L$ be a finite lattice such that $\mathcal{A}(L)$ is independent. Then, the M\"obius function $M$ of $L$ is
$ M(x)= (-1)^{\card B}$ if $x:=\bigvee B$, for some $B \subseteq \mathcal{A}(L)$, and $ M(x)=0$ otherwise.
\end{theo}
As a corollory, for all $\vec{v}:=(v_1, \ldots, v_n) \in \N^n$, $M(\vec{v}):=(-1)^{t_{\vec{v}}}$, if and only if for all $i$, $v_i:=1$ or $i$, and where $t_{\vec{v}}$ is the number of coordinates such that $v_i:=i \not=1$. Else, $M(\vec{v}):=0$. 
\subsubsection{Realisation of the grafting operation}
In the sequel, for all $n>0$, $\vec{n}:=(1,2,3,4, \ldots,n)$, $\vec{0}:=(0)$ and $\vec{1_n}:=(1,1,1, \ldots,1)$.
Fix $n,m \not=0$ and $\vec{v} \in \N^n$ and 
$\vec{w} \in \N^m$. The grafting operation is a map,
$$ \vee: \N^n \times \N^m \xrightarrow{} \N^{m+n+1}, \ \ \  (\vec{v}, \vec{w}) \mapsto \vec{v} \vee \vec{w}:=
(\vec{v},1, w_1+n+1, \ldots, w_m +n+1):=(\vec{v},1,v+1+\vec{w}),$$
where for all $k \in \Na$ and $\vec{w} \in \N^m$, $m>0$, the notation $k+\vec{w}$ stands for $(w_1+k, \ldots, w_m +k)$,  $k+\vec{0}:=\vec{0}$ and where by abuse of notation $v$ denotes the number of coordinates in $\vec{v}$, (\textit{i.e.}, $n$ in this case).
In the sequel, we give the name $(0)$ to the tree $\treeO$ and 
$(1)$ to the tree $Y:=\treeA$. Hence, $M((1)):=+1$. 
By convention, if $\vec{v}:=(v_1, \ldots, v_n) \in \N^n$ and $n \not=0$, then $\vec{v} \vee (0):=(\vec{v},1)$, $(0) \vee \vec{v} :=(1, 1+\vec{v})$ and $(0) \vee (0):= (1)$.
We extend the M\"obius function to $(0)$ by setting:
$ M(\vec{v} \vee (0)):=M(\vec{v},1):=M(\vec{v}),$ and $M((0) \vee \vec{v}):=M((1, 1+ \vec{v}):=0,$ unless $\vec{v}:=\vec{n},$
for all $\vec{v} \in \N^n$ and $n>0$.
\begin{prop}
\label{ineg}
Let $\pi \in Y_n$ and $\tau \in Y_m$. Then, $\vec{name}(\pi \vee \tau) = \vec{name}(\pi) \vee \vec{name}(\tau)$. (The map $\vec{name}$ is a grafting morphism.) Moreover, $M(\vec{name}(\pi \vee \tau)) :=(-1)^m M(\vec{name}(\pi))$ if $\tau:=\vec{m}$ and zero otherwise.
Moreover,
if $\vec{v},\vec{w} \in (\N^n,<)$ and $\vec{u}, \vec{z} \in (\N^m,<)$, then
$ \vec{v} \leq \vec{w}, \ \vec{u} < \vec{z} \ \ (\textrm{or} \ \vec{v} < \vec{w}, \ \vec{u} \leq \vec{z}) \Leftrightarrow \vec{v} \vee \vec{u} < \vec{w} \vee \vec{z}.$
\end{prop}
\Proof
Let $\pi \in Y_n$ and $\tau \in Y_m$. The tree $\pi$ gives a unique complete expression,
$Exp(\pi):=(^{p_1}x_1(^{p_2}x_2 \ldots (^{p_{n}}x_{n} x_{n+1},$ 
where $p_i$ is the number of $i$ in the name of $\pi$. Similarly
for  $\tau$, set  $Exp(\tau):=(^{p'_1}x_1(^{p'_2}x_2 \ldots (^{p'_{m}}x_{m} x_{m+1}$. Their grafting  gives $$((^{p_1}x_1(^{p_2}x_2 \ldots (^{p_{n}}x_{n} x_{n+1}(^{p'_1}x_1(^{p'_2}x_2 \ldots (^{p'_{m}}x_{m} x_{m+1},$$
which once renamed in a complete expression of $\bra \ x_1, \ldots, x_{n+1},x_{n+2}, \ldots, x_{n+1+m+1}, ( \ \ket$ gives
$((^{p_1}x_1(^{p_2}x_2 \ldots (^{p_{n}}x_{n} x_{n+1}(^{p'_1}x_{n+1+1}(^{p'_2}x_{n+1+2} \ldots (^{p'_{m}}x_{n+1+m} x_{n+1+m+1}.$
Observe that $\vec{name}(\pi \vee \tau)_{n+1}:=1$, giving the first claim. For computing the M\"obius function, observe that if $1 \leq w_i < i$, then $n+2 \leq \vec{name}(\pi \vee \tau)_{i+n+1} < i+n+1$. Without forgetting $w_1:=1$, which
becomes $\vec{name}(\pi \vee \tau)_{1+n+1}:=n+2>1$, we obtain $M(\vec{name}(\pi \vee \vec{m})) :=(-1)^{m}M(\vec{name}(\pi))$ and zero otherwise. The last claim is straightforward.
\eproof

\noindent
The localisation of $v_i:=1, \ i \not=1,$ inside a given name of a binary tree reveals grafting operations since the most right parenthesis at the right hand side of $x_i$ closes a left parenthesis open at $x_1$ giving thus a complete expression in $\bra \ x_1, \ldots, x_n, ( ,)\ \ket$. There are also `hidden' graftings 
due to the translation of $n+1$ in the vector $\vec{w}$.
There exists a trivial partial order in $\Na^n$ by declaring that 
$\vec{v} \leq \vec{w} \Leftrightarrow \ \forall \ 1\leq i \leq n, \ v_i \leq w_i,$ inducing so a trivial partial order on $\N^n$. As already mentioned, there exists a partial order on binary trees, often called the Tamari order, induced by the relation $(\tau_1 \vee \tau_2) \vee \tau_3 \leq \tau_1 \vee (\tau_2 \vee \tau_3),$ for any trees $\tau_1, \tau_2, \tau_3$. 
Equip $Y_n$ with the Tamari order. Then, for all 
$\pi,\tau \in Y_n$, 
$ \pi < \tau$ if and only if $\vec{name}(\pi) < \vec{name}(\tau).$
There is on the symmetric group $S_n$, a partial order $\leq_{Bruhat}$ called the weak-Bruhat order. From \cite{LRbruhat}, there is a surjective map, $s \mapsto Y_s$, mapping permutations of $S_n$ to $n$-trees of $Y_n$. Equipped with the Tamari order, it is proved that 
$s \leq_{Bruhat} s' \Leftrightarrow Y_s \leq Y_{s´}$. Therefore, the  weak-Bruhat order of the symmetric group $S_n$ is nothing else that the trivial partial order on $\N^n$. As $S_n$ is in bijection with $\Nb^n$, it might be interesting to find an order preserving code between permutations and vectors of $\Nb^n$.
The Tamari order is represented for $(Y_2,<), (\N^2,<)$,
$ \treeAB \ (1,1) \rightarrow \treeBA \ (1,2)$ and
for $(Y_3, <)$ or $(\N^3,<)$,
\begin{small}
\begin{center}
$
\begin{array}{ccccc}
 & & \treeABC (1,1,1) & &\\
 & \swarrow & &\searrow &  \\
& & & &\treeBAC (1,2,1) \\
 \treeACA (1,1,3) & & & & \downarrow \\
 & & & & \treeCAB (1,2,2)\\
 & \searrow & &\swarrow &  \\
 & & \treeCAB (1,2,3) & & \\
\end{array}
$
\end{center}
\end{small}
For all $n>0$, we mention the existence of convex polytopes, so-called
Stasheff polytopes or associahedrons \cite{Stasheff, Lodayst}, denoted by $\mathcal{K}^{n}$ and whose vertices are indexed by the binary trees of  $Y_{n+1}$. (Just above $\mathcal{K}^{1}$ and $\mathcal{K}^{2}$ are represented.) The vector formulation give a majoration of the number of paths between two vertices. Indeed,
if $\pi,\tau \in (Y_n,<)$ with $\pi < \tau$, then, the number of paths from $\pi$ to $\tau$ in $\mathcal{K}^{n-1}$ is less or equal to $\Pi_{i=1} ^n \ (name(\pi)_i -name(\tau)_i +1)$. Another way to check if a vector of $\Nb^n$ is the name of a tree is the following. Fix 
$\vec{v} \in \Nb^n$ and take the highest coordinate such that $v_i:=1$. We get a unique decomposition $\vec{v}:=(\vec{v}_l,1, v_l+1+ \vec{v}_r)$. The vector
$\vec{v}$ is the name of a tree if and only if so are $\vec{v}_l$
and $((-v_l -1) + \vec{v}_r)$.
The grafting operation can be extended by bilinearity to
$K\N^{\infty}:= \bigoplus_{n \geq 0} K\N^n$. In the sequel,
we set $K\N^{\infty}_*:= \bigoplus_{n \geq 1} K\N^n$, $ \N^{\bullet}:= \cup_{n \geq 0} \N^n$ and $ \N^{\bullet}_*:= \cup_{n \geq 1} \N^n$. 
\subsubsection{Coding the over and under operations}
\label{over}
Before going on, recall that an associative \textit{$L$-algebra} is a $K$-vector space $A$ equipped with two binary operations $\nearrow, \ \nwarrow: A^{ \otimes 2} \xrightarrow{} A$ and obeying three constraints. The two operations are associative and verify the `link': $(x\nearrow y)\nwarrow z:= x\nearrow (y\nwarrow z)$. From a coalgebraic point of view, $L$-coalgebras have been introduced on graphs in \cite{Coa, codialg1, dipt}.
In \cite{LRbruhat}, J.-L. Loday and M. Ronco introduced the operations \textit{over} and \textit{under} on trees, denoted respectively by
$\nearrow, \nwarrow: Y_n \times Y_m \xrightarrow{} Y_{n+m}$, for all $n,m \not=0$, where $\pi \nearrow \tau$ is the tree $\tau$ with its most left leaf  identified with the root of $\pi$ and where $\pi \nwarrow \tau$ is the tree $\pi$ with its most right leaf  identified with the root of $\tau$. These two operations have a common unit which is $\treeO$.
To define the analogue of these two operations on vectors,  consider the map $\triangleright: \Na \times \N^n \xrightarrow{} \N^n,$ $k \triangleright \vec{v}:= (k\tilde{+}v_1, \ldots k \tilde{+} v_n)$, where $k \tilde{+} v_i:= k+v_i$, for $v_i \not=1$ and $k \tilde{+} 1:=1$, (otherwise stated, 1 is a right anihilator for  the operation $\tilde{+}$).
\begin{prop}
\label{ovund}
Fix $n,m \not=0$ and $\vec{v} \in \N^n$ and $\vec{w} \in \N^m$. The binary operations $\nearrow, \nwarrow: \N^n \times \N^m \xrightarrow{} \N^{n+m}$ defined as follows:
$ \vec{v} \nearrow \vec{w}:= (\vec{v},  v \triangleright \vec{w}), \ \ \textrm{and} \ \ \vec{v} \nwarrow \vec{w}:=(\vec{v}, v+\vec{w}),$
turn $\N^\bullet_*$, (resp. $K\N^{\infty}_*$) into an associative $L$-monoid (resp. an associative $L$-algebra).  The map $\vec{name}$ is a morphism of associative $L$-monoids, (resp. of associative $L$-algebras). Moreover, the M\"obius
function has a simple expression:
$M(\vec{v} \nearrow \vec{w}):= M(\vec{v}) M(\vec{w}) \ \ and \ \
M(\vec{v} \nwarrow \vec{w}):= (-1)^m M(\vec{v}), \ \textrm{if} \ \vec{w}:=\vec{m}, \ \textrm{else} \ 0.$
\end{prop}
\Proof
Fix $n,m \not=0$ and $\vec{v} \in \N^n$ and $\vec{w} \in \N^m$. Their complete expression 
$Exp(\vec{v})$ (resp. $Exp(\vec{w})$) is of the form $(^{p_1}x_1(^{p_2}x_2 \ldots
(^{p_n}x_nx_{n+1}$, (resp. $(^{p'_1}x_1(^{p'_2}x_2 \ldots
(^{p'_m}x_mx_{m+1}$). The associated trees are $Tree(\vec{v})$ and $Tree(\vec{w})$. However, $Tree(\vec{v}) \nearrow Tree(\vec{w})$ has the expression,
$(^{p'_1}Exp(\vec{v})(^{p'_2}x_2 \ldots
(^{p'_m}x_mx_{m+1}.$
Observe that  $(\vec{v} \nearrow \vec{w})_{n+1}$ corresponds to $x_{n+1}$ thus is equal to 1. Observe also that the left parentheses of $Exp(\vec{w})$ do not move during this operation. We have to take into account the shift of the coordinate $j$ of $\vec{w}$ of an amount of $v$ --corresponding to the degree of the tree $Tree(\vec{w})$-- for all $w_j \not=1$. For
$w_j =1$, the most right parenthesis at the right hand side of $x_j$ still close a left parenthesis at the left hand side of $x_1 \equiv Exp(\vec{v})$. Therefore, for those $j$, $(\vec{v} \nearrow \vec{w})_{j}:=1.$
This gives the vector $(\vec{v},1,v \tilde{+} w_2, v \tilde{+} w_3, \ldots, v \tilde{+} w_m):=(\vec{v}, v \triangleright \vec{w}).$ 
The second operation is easier since all the $w_j$ have to be shifted by $n:=v$. We extend easily these two operations to $\N^\bullet_*$ and to $K\N^{\infty}_*$ (by bilinearity for the second case). Observe then, that $\nearrow$ and $\nwarrow$ are associative and the equality, $\vec{u} \nearrow(\vec{v} \nwarrow \vec{w})=(\vec{u} \nearrow \vec{v} )\nwarrow \vec{w},$ holds, giving an associative $L$-monoidal structure
to $\N^\bullet_*$ or an associative $L$-algebra structure to $K\N^{\infty}_*$. Extend the map $\vec{name}$ by linearity, we obtain an isomorphism of associative $L$-algebras
between  $KY^{\infty}_*$ and $K\N^{\infty}_*$. Concerning
the identities on the M\"obius function, observe that 
$M(\vec{v} \nearrow \vec{w}):= M(\vec{v}) M(\vec{w}) $ since
$w_i:=i \not=1$ if and only if $(\vec{v} \nearrow \vec{w})_ {i+v}:= i+v$ and $v_i:=1$ becomes $(\vec{v} \nearrow \vec{w})_ {i+v}:=1$. However $M(\vec{v} \nwarrow \vec{m}):= (-1)^{m}M(\vec{v}) $ since we have also to take into account $w_1:=1 $ becoming $(\vec{v} \nearrow \vec{w})_ {1+v}:= 1+v$ and for $\vec{w} \not= \vec{m}$, $M(\vec{v} \nwarrow \vec{w})=0$.
\eproof
\begin{coro}
\label{zzt}
For all $\vec{u}_1, \vec{u}_2  \in (\N^n,<), \ \vec{v}_1, \vec{v}_2 \in (\N^m, <)$ and $\vec{w}_1, \vec{w}_2 \in (\N^p, <)$,
$\vec{u}_1\leq \vec{u}_2,  \ \vec{w}_1 \leq \vec{w}_2, \ \vec{v}_1 \leq \vec{v}_2 \Leftrightarrow \vec{u}_1 \nearrow \vec{v}_1 \nwarrow \vec{w}_1 \leq \vec{u}_2 \nearrow \vec{v}_2 \nwarrow \vec{w}_2,$
where the presence of a strict inequality on the left hand side induces a strict one in the right hand side.
Moreover,
$\vec{v} \nearrow \vec{w} \leq \vec{v} \nwarrow \vec{w}$, holds.
\end{coro}
\Proof
The proof is complete by using Proposition \ref{ovund}.
\eproof

Since $\vec{v}:=\vec{v}_l \vee \vec{v}_r= \vec{v}_l \nearrow (1) \nwarrow \vec{v}_r$, it is straightforward to prove that the free $L$-algebra over one generator $x$ is isomorphic to $(KN^{\infty}_*, \nearrow, \nwarrow)$ by mapping $x$ to the generator $(1)$ (see also \cite{Pir}). 
By anticipating the ideas of J.-L. Loday (explained in details below), one can convert operations $\nearrow, \ \nwarrow$ into set operations called $L$-additions denoted by $+_{_{\nearrow}}, +_{_{\nwarrow}}: \N^{n} \times  \N^{m} \xrightarrow{} \N^{n+m}$ where $\vec{v} +_{_{\nwarrow}} \vec{u}:= \vec{v} \nwarrow \vec{u}$ and $\vec{v} +_{_{\nearrow}} \vec{u}:= \vec{v} \nearrow \vec{u}$. These additions are associative and noncommutative. Similarly, there is a notion of $L$-multiplication. As $(KN^{\infty}_*, \nearrow, \nwarrow)$ is the free $L$-algebra on the generator $(1)$, one can uniquely write
any name of binary trees \textit{via} the operations $\nearrow$ and $\nwarrow$ and $(1)$. Such a formula for a vector $\vec{v}$, is called its universal expression and is denoted by $\varpi_{\vec{v}}((1))$, obtained by the following induction $\varpi_{\vec{v}}((1)):=\varpi_{\vec{v}_l}((1))\nearrow (1) \nwarrow \varpi_{\vec{v}_r}((1))$. For instance, $(1,1,3):=(1)\nearrow (1) \nwarrow (1):=\varpi_{(1,1,3)}((1)).$ The $L$-multiplication of $\vec{u} \in \N^{n}$ by $\vec{v} \in \N^{m}$ is by definition: $\vec{u} \tilde{\ltimes} \vec{v}:= \varpi_{\vec{u}}(\vec{v}) \in \N^{nm}$. For instance: $(1,1,3)\tilde{\ltimes} \vec{v}:=(\vec{v})\nearrow (\vec{v}) \nwarrow (\vec{v})$. Therefore, any name of  $\N^{m}$, where $m$ is a prime number, will be prime for the
$L$-arithmetics. Consider now the $K$-vector space $K[X]_L$ spanned by
$\{X^{\vec{v}}, \vec{v} \in (N^{\bullet}_*, +_{_{\nearrow}}, +_{_{\nwarrow}}, \tilde{\ltimes})\}$. This is the free $L$-algebra over the generator $X^{(1)}$ where as expected, operations are defined by $X^{\vec{u}} \nearrow X^{\vec{v}}:= X^{\vec{u} +_{_{\nearrow}} \vec{v}}, \ X^{\vec{u}} \nwarrow X^{\vec{v}}:= X^{\vec{u} +_{_{\nwarrow}} \vec{v}}$ and $(X^{\vec{u}})^{\vec{v}}:=X^{\vec{u}\tilde{\ltimes} \vec{v}}$, imitating the usual ploynomial algebra on one variable endowed with the usual arithmetics over $\Na$. There is also a dendriform involution $\dagger$, described in Subsubsection \ref{dendinvosubsect}. We summarize our investigation by the following theorem.
\begin{theo}
The set $N^{\bullet}_*$ equipped with the $L$-additions, $+_{_{\nearrow}}$ and $+_{_{\nwarrow}}$ with the $L$-multiplication $\tilde{\ltimes}$ and with the dendriform involution $\dagger$ is an involutive graded $L$-monoid. The $L$-multiplication is left distributive, associative though noncommutative. For any names of trees, $\vec{u}, \vec{v}$, 
$(\vec{u} +_{_{\nearrow}} \vec{v})^{\dagger}:= \vec{v}^{\dagger} +_{_{\nwarrow}} \vec{u}^{\dagger}$,
$(\vec{u} +_{_{\nwarrow}} \vec{v})^{\dagger}:= \vec{v}^{\dagger} +_{_{\nearrow}} \vec{u}^{\dagger}$ and $(\vec{u} \tilde{\ltimes} \vec{v})^{\dagger}:= \vec{u}^{\dagger} \tilde{\ltimes} \vec{v}^{\dagger}$, hold.
Moreover, equipped with the dendriform involution, the $K$-vector space
$K[X]_L$ spanned by
$\{X^{\vec{v}}, \vec{v} \in (N^{\bullet}_*, +_{_{\nearrow}}, +_{_{\nwarrow}}, \tilde{\ltimes})\}$ is the free involutive associative $L$-algebra over the generator $X^{(1)}$.
\end{theo}
\begin{prop}
\label{dd}
Fix $n,m \not=0$ and $\vec{x}, \vec{y} \in (N^{n},<)$ and $\vec{a}, \vec{b} \in (N^{m},<)$.
With regards to the trivial partial order, 
the map $\varpi_{\vec{x}}: N^{m} \xrightarrow{} N^{nm}$, is a lattice morphism, \textit{i.e.,} $\vec{x} \tilde{\ltimes} \vec{a} < \vec{x} \tilde{\ltimes}  \vec{b} \Leftrightarrow \vec{a} < \vec{b}$ and
the map $ \tilde{\ltimes} \vec{a}: N^{n} \xrightarrow{} N^{nm}$ is also a lattice morphism, \textit{i.e.,} $\vec{x} \tilde{\ltimes} \vec{a} < \vec{y} \tilde{\ltimes}  \vec{a} \Leftrightarrow \vec{x} < \vec{y}$.
\end{prop}
\Proof
Keep notation of Proposition \ref{dd}.
For the first claim, proceed by induction. It is true for $\vec{x}:=(1)$, for $\vec{x}:=(1,1)$ and for  $\vec{x}:=(1,2)$.  By Proposition \ref{zzt}, 
$\vec{x} \tilde{\ltimes} \vec{a} < \vec{x} \tilde{\ltimes} \vec{b} \Leftrightarrow \vec{x}_l \tilde{\ltimes} \vec{a} +_{_{\nearrow}} a +_{_{\nwarrow}} \vec{x}_r \tilde{\ltimes} \vec{a}
< \vec{x}_l \tilde{\ltimes} \vec{b} +_{_{\nearrow}} \vec{b} +_{_{\nwarrow}} \vec{x}_r \tilde{\ltimes} \vec{b} \Leftrightarrow \vec{a} < \vec{b}, \ \ \vec{x}_l \tilde{\ltimes} \vec{a} < \vec{x}_l \tilde{\ltimes} \vec{b}, \ \ \vec{x}_r \tilde{\ltimes} \vec{a} <  \vec{x}_r \tilde{\ltimes} \vec{b}.
$ 
The proof is complete by induction. Concerning the second claim, if $\vec{x} < \vec{y}$, then there exist say $k$ Tamari moves between the trees associated with $\vec{x}$ and $\vec{y}$. Suppose $k=1$. Then, in the definitions of $\vec{x}, \vec{y}$, this means the existence of three vectors say $\vec{v}_1,\vec{v}_2,\vec{v}_3  $ such as we have $\ldots (\vec{v}_1\vee \vec{v}_2) \vee \vec{v}_3 \ldots < \ldots \vec{v}_1\vee (\vec{v}_2 \vee \vec{v}_3) \ldots$. Therefore, we obtain,
$\ldots (\vec{v}_1\nearrow (1) \nwarrow \vec{v}_2) \nearrow (1) \nwarrow \vec{v}_3 \ldots < \ldots \vec{v}_1 \nearrow (1) \nwarrow (\vec{v}_2 \nearrow (1) \nwarrow \vec{v}_3) \ldots $. The second claim holds for $k=1$, and thus for all $k$.
\eproof

\noindent
As binary trees are considered up to isotopies, the operations $\nearrow, \  \nwarrow$ have a common unit which is $(0) \equiv \treeO$. However,
 the link axiom of $L$-algebras is not compatible with this unit since it forces $\nearrow \ = \ \nwarrow$. Using 
the trivial partial order, we will exhibit an associative operation, sum of two nonassociative operations obeying three axioms. This operation has first been introduced by J.-L. Loday and M. Ronco by using technics in permutation groups \cite{LRbruhat}. One of the main advantages of our coding is to give easier proofs to these results.
\begin{prop}
\label{dend}
The following binary operation,
$ \star: K\N^{\infty}_* \otimes K\N^{\infty}_* \longrightarrow  K\N^{\infty}_*, \ \ \vec{v} \otimes  \vec{w} \mapsto
\vec{v} \star \vec{w}:=\sum_{\vec{v}  \nearrow \vec{w} \leq \vec{t} \leq \vec{v}  \nwarrow \vec{w}} \ \vec{t},$
is associative. Moreover, $ \vec{u} \star \vec{v} \star \vec{w}:=\sum_{\vec{u} \nearrow \vec{v}  \nearrow \vec{w} \leq \vec{t} \leq \vec{u} \nwarrow \vec{v}  \nwarrow \vec{w}} \ \vec{t}$ 
and $\vec{v} \star (0)=\vec{v}=(0) \star \vec{v}$ hold for all
$\vec{u}, \vec{v}, \vec{w} \in K\N^{\infty}_*$.
\end{prop}
\Proof
Let $\vec{u} \in \N^{p}, \vec{v} \in \N^{n}$ and $\vec{w} \in \N^{m}$, with $p,n,m \not=0$. Write down $\vec{u} \star (\vec{v} \star \vec{w})$ to obtain the square $S_1$ --here the dots mean $< \ldots <$--:

$
\begin{array}{ccc}
(\vec{u} \nearrow \vec{v}) \nearrow \vec{w} & \ldots & (\vec{u} \nearrow \vec{v}) \nwarrow \vec{w}\\ 
 \vdots & \vdots &\vdots \\
(\vec{u} \nwarrow \vec{v}) \nearrow \vec{w} & \ldots & (\vec{u} \nwarrow \vec{v}) \nwarrow \vec{w}.
\end{array}
$

\noindent
Write down $(\vec{u} \star \vec{v}) \star \vec{w}$ to obtain the square $S_2$:

$
\begin{array}{ccc}
\vec{u} \nearrow (\vec{v} \nearrow \vec{w}) & \ldots & \vec{u} \nearrow (\vec{v} \nwarrow \vec{w})\\ 
 \vdots & \vdots &\vdots \\
\vec{u} \nwarrow (\vec{v} \nearrow \vec{w}
) & \ldots & \vec{u} \nwarrow (\vec{v} \nwarrow \vec{w}).
\end{array}
$

\noindent
Use associativity of $\nearrow$ and $\nwarrow$ and
$(\vec{u} \nwarrow \vec{v}) \nearrow \vec{w} := (\vec{u}, u + \vec{v}, (u +  v) \triangleright  \vec{w})) <
(\vec{u}, u + \vec{v}, u + ( v \triangleright \vec{w}))=:\vec{u} \nwarrow (\vec{v} \nearrow \vec{w}
)$ to complete the proof. The last claim is obvious since $(0)$ is by definition a unit for the operations $\nearrow$ and $\nwarrow$.
\eproof

\noindent
The sum in the definition of the associative product $\star$ can be split into two parts corresponding to two operations. 
\begin{prop}
\label{splii}
Let $\vec{v},  \vec{w}$ be names of some trees. Then, the set $I:=\{\vec{u}; \ \vec{v} \nearrow  \vec{w} \leq \vec{u} \leq    \vec{v} \nwarrow  \vec{w} \}$ splits into two disjoint subsets: $I_1:=\{\vec{u}; \ (\vec{v}, v_l +1+v_r \triangleright  \vec{w}) \leq \vec{u} \leq    (\vec{v}, v+  \vec{w} )   \}$ and $I_2:=\{\vec{u}; \ (\vec{v}, v \triangleright  \vec{w}) \leq \vec{u} \leq   ( \vec{v},   v+ \vec{w}_l,1,  v+1+ w_l+ \vec{w}_r    \}$.
\end{prop}
\Proof
First of all, observe that $(\vec{v}, v \triangleright  \vec{w}) :=
(\vec{v}, v \triangleright  \vec{w}_l,1, v+1+w_l+\vec{w}_r)$. Therefore, we have only to compare $(v \triangleright  \vec{w}_l,1)$ and $(v +  \vec{w}_l,1)$ in $I_2$.
Similarly, concerning $I_1$, observe that $(\vec{v}, v_l +1+v_r \triangleright  \vec{w}) := (\vec{v}, v_l +1+v_r \triangleright  \vec{w}_l,v_l +2, v_l +1+v_r \triangleright  [\vec{w}_r + w_l +1]) = (\vec{v}, v_l +1+v_r \triangleright  \vec{w}_l,v_l +2, v+1+ w_l +\vec{w}_r )$ and  $(\vec{v}, v+  \vec{w} ) :=(\vec{v}, v+  \vec{w}_l, v+1, v+ w_l +1+ \vec{w}_r ) $. Therefore, we have to compare the vector $(v_l +1+v_r \triangleright  \vec{w}_l,v_l +2)$ with $(v+  \vec{w}_l, v+1)$.
As $\vec{v}_l$ represents a complete expression, jumps of coordinates situated after $v_l$ cannot take values below $v_l$. From this remark, one obtains that $I_1$ and $I_2$ are disjoint and $I_1 \cup I_2=I$.
\eproof

\noindent
We recover the dendriform dialgebra introduced in \cite{Loday} from a vectorial framework. Recall that
a $K$-vector space $E$ is a {\it{dendriform dialgebra}} \cite{Loday} if it is equipped with 2 binary operations $\prec$ and $\succ$ satisfying the following axioms
for all $x,y \in E$:
$$(x \prec y )\prec z = x \prec(y \star z), \ \ \
(x \succ y )\prec z = x \succ(y \prec z), \ \ \
(x \star y )\succ z = x \succ(y \succ z), \ \ \ $$
where, by definition, $x \star y :=x  \prec y +x \succ y$, for all $x,y \in E$, where $\star$ turns out to be associative. We now give a different proof
of the following theorem appearing in \cite{Loday} and \cite{LRbruhat}.
\begin{theo}
Equip $K\N^{\infty}_*$ with two binary operations $\prec$ and $\succ$, defined as follows:
$ \vec{v} \prec \vec{w}:=\vec{v}_l \vee (\vec{v}_r \star \vec{w})$ and $ \vec{v} \succ \vec{w}:=(\vec{v} \star \vec{w}_l) \vee \vec{w}_r ,  \ \forall \ \vec{v},\vec{w} \not=(0).$
Then, $(K\N^{\infty}_*, \prec, \succ)$ is a dendriform dialgebra generated by $(1)$. This space can be augmented by requiring $ \vec{v} \prec (0):= \vec{v} =: (0) \succ \vec{v} $ and
$ \vec{v} \succ (0):= 0=: (0) \prec \vec{v}$, for $\vec{v} \not=(0)$.
Equipped with these operations, $K\N^{\infty}$ is still a dendriform dialgebra with $ \vec{v} * (0)=(0)*\vec{v}=\vec{v}$, for all $\vec{v} \in K\N^{\infty}$.
\end{theo}
\Proof
Observe that $(0)*(0):=(0)$ but that $(0) \prec (0)$ and $(0) \succ (0)$ are not defined. The associative operation $\star$ of Proposition \ref{dend}
is associative and splits into two operations $\prec, \succ$ according
to Proposition \ref{splii} and defined in this theorem.
Let us prove the first axiom of dendriform dialgebras. Fix $\vec{u}, \vec{v}, \vec{w} \not=0$. Then,
$$ (\vec{u} \prec \vec{v })  \prec \vec{w} := (\vec{u}_l \vee (\vec{u}_r \star \vec{v }))  \prec \vec{w} :=\vec{u}_l \vee (\vec{u}_r \star \vec{v }\star \vec{w} ):= \vec{u}_l \vee (\vec{u}_r \star (\vec{v }\star \vec{w}) ) := \vec{u} \prec (\vec{v } \star \vec{w}). $$
By induction, the vector $(1)$ is the generator of $(K\N^{\infty}_*, \prec, \succ)$ since,
\begin{eqnarray*}
\vec{v}_1 \vee \vec{v}_2  &:=& (1)   \ \ \ if  \ \vec{v}_1 =(0) = \vec{v}_2, \\
&:=& (1) \prec \vec{v}_2 \ \ \ if \  \vec{v}_1 = (0) \not= \vec{v}_2, \\
&:=& \vec{v}_1 \succ (1)  \ \ \ if \  \vec{v}_1 \not= (0) = \vec{v}_2, \\
&:=& \vec{v}_1 \succ (1) \prec \vec{v}_2  \ \ \ if \  \vec{v}_1 \not= (0) \not= \vec{v}_2.
\end{eqnarray*}
The second claim is obtained by checking that axioms of dendriform dialgebras are compatible with the unit action defined in this theorem.
\eproof
\begin{theo}[Loday \cite{Loday}]
The  $K-$vector space $(K\N^{\infty}_*,\prec,\succ)$ is the free dendriform dialgebra
on the generator $(1)$.
\end{theo}
\Proof
This a reformulation of a result in \cite{Loday}.
\eproof

\noindent
As a corollary, there exists a \textit{universal expression}, denoted by $\omega_{\vec{v}}(1)$, of $\vec{v} \in \N^n$
as a composition of $n$ copies of (1) with $\prec$ and $\succ$. Set $\omega_{(0)}(1):=0$ and of course $\omega_{\vec{v}}(1):=\omega_{\vec{v}_l}(1)\succ (1) \prec \omega_{\vec{v}_r}(1) $. For instance, $\omega_{(121)}(1):=((1) \prec (1))\succ (1).$
So defined, $(K\N^\infty_*:=\otimes_{n>0} K\N^n, \ \prec, \succ, \nearrow, \nwarrow)$
is another representation (more tractable) of the free dendriform dialgebra on one generator $(1)$, equip with an extra-structure of associative $L$-algebra  $(K\N^\infty, \nearrow, \nwarrow)$ whose basis encodes binary trees in a compatible way with
the Tamari order underlying the definitions of the operations
$\prec$ and $\succ$. One of the main advantages of this coding lies in a slight reformulation of
arithmetree in terms of vectors. 
\subsection{Recall of arithmetree on planar binary trees}
\label{arithmdendsec}
After these slight reformulations of the constructions developed in \cite{Loday, LRbruhat}, let us recall a deep
notion introduced by J.-L. Loday. We follow  \cite{Lodayarithm}. A {\it{grove}} is simply a non-empty subset of $Y_n$, \textit{i.e.}, a disjoint union of binary trees with same degree such that each tree appears only once. The set of groves over $Y_n$ is denoted by $\YY_n$ and is of cardinal $2^{c_n}-1$. For instance in low degrees,
$$ \YY_0:=\{ \treeO \}, \YY_1:=\{ \treeA \},  \YY_2:=\{ \treeAB,
\treeBA,  \treeAB \cup \treeBA \}.$$
Similarly, we define $\Ng^n$ in the same way. Instead of  binary trees, we work with a more tractable set $\Ng^n$, which are the names of groves of $\YY_n$. Hence,
$ \Ng^0:=\{(0) \}, \Ng^1:=\{ (1) \},  \Ng^2:=\{ (1,1),
(1,2), (1,1)  \cup (1,2) \},$
and continue to call grove such a union of vectors.
The idea is to convert the associative operation $\star$ in Proposition \ref{dend} into an addition with values in groves.
\subsubsection{The dendriform addition}
\begin{defi}{[Dendriform addition \cite{Lodayarithm}]}
The \textit{dendriform addition} of two vectors $\vec{v}$ and $\vec{w}$ associated with some planar binary trees is defined by:
$$ \vec{v} \dotplus \vec{w} := \bigcup_{\vec{v} \nearrow \vec{w} \leq \vec{u} \leq \vec{v} \nwarrow \vec{w}}     \vec{u}. $$
This is extended to groves by distributivity of both sides, \textit{i.e.},
$ \cup_i \vec{v}_i \dotplus \cup_j \vec{w}_j:=  \cup_{ij} (\vec{v}_i  \dotplus  \vec{w}_j),$ which has a meaning thanks to Theorem \ref{sum}.
For instance: $(1) \dotplus (1):= (1,1) \cup (1,2)$ or at the level of binary trees $\treeA  \dotplus \treeA = \treeAB  \cup \treeBA$.
\end{defi}
\textbf{Warning:}
Though associative, the dendriform addition is not commutative.  
$(1) \dotplus (1,1):= (1,1,1) \cup  (1,2,1) \cup (1,2,2)$ giving a grove different from $(1,1) \dotplus (1):= (1,1,1) \cup (1,3,1)$.
\begin{lemm}
\label{encadr}
Let $\vec{w} \in \N^{n+m}$. Then, there exists unique  $\vec{u} \in \N^{n}$ and $\vec{v}  \in \N^{m}$ such that:
$\vec{u} \nearrow \vec{v} \leq \vec{w} \leq \vec{u} \nwarrow \vec{v}.$
\end{lemm}
\Proof
(Compare to \cite{Lodayarithm}, Prop. 2.3. and Corol. 2.4.). Recall that for $\vec{u} \in \N^{n}$ and $\vec{v}  \in \N^{m}$, we get:
$\vec{u} \nearrow \vec{v}:= (\vec{u}, u \triangleright \vec{v})$ and $\vec{u} \nwarrow \vec{v}:= (\vec{u}, u+\vec{v})$. Take the first $n$ coordinates of $\vec{w} \in \N^{n+m}$. This gives
a unique vector $\vec{u} \in \N^{n}$ according to Proposition \ref{err}. Consider the vector $\vec{v_1}$ defined by $\vec{v_1}:=(w_{n+1}, \ldots, w_{n+m})$. Make the translation of $-n:=-u$ to obtain $\vec{v_1}-u:=(w_{n+1}-u, \ldots, w_{n+m}-u)$. The vector $\vec{v}  \in \N^{m}$ we are looking for is obtained by replacing all negative or null coordinates by $1$. Observe that:
$\vec{u} \nearrow \vec{v}:= (\vec{u}, u \triangleright \vec{v}) \leq \vec{w} \leq (\vec{u}, u+\vec{v}):=\vec{u} \nwarrow \vec{v}.$
\eproof

\noindent
We now simplify the proof of the following theorem.
\begin{theo}[Loday, \cite{Lodayarithm}]
 \label{sum}
The dendriform addition of two groves is still a grove, \textit{i.e.,}
$  \dotplus: \Ng^n \times \Ng^m \xrightarrow{} \Ng^{n+m}.$
\end{theo}
\Proof
A priori, it is not immediate that binary trees appearing in the union defining the dendriform addition are all different. Nevertheless, consider the total grove $\u{n+1}:=\cup_{\vec{w} \in \N^{n+1}} \ \vec{v}$. By applying Lemma \ref{encadr}, observe that $\u{n + 1}:=\cup_{\vec{v} \in \N^{n}} \cup_{\vec{v} \nearrow (1)\leq \vec{w} \leq \vec{v} \nwarrow (1)} \ \vec{w} := \u{n} \dotplus \u{1}$. Apply associativity of the dendrifrom addition and induction to obtain $\u{n} \dotplus \u{m}:=\u{n} \dotplus \u{1}  \dotplus \u{1} \ldots  \dotplus \u{1}:=\u{n+m}. $
\eproof
\begin{prop}[Left and right cancellations]
Let $\vec{u},\vec{v} \in \N^{m},\vec{w} \in \N^{n}$. Then,
$\vec{u} \dotplus \vec{v} = \vec{u} \dotplus \vec{w} \Leftrightarrow \vec{v}= \vec{w}$ and
$\vec{v} \dotplus \vec{u} = \vec{w} \dotplus \vec{u} \Leftrightarrow \vec{v}= \vec{w}$.
\end{prop}
\Proof
Apply Proposition \ref{ovund} to conclude.
\eproof

\noindent
\textbf{[Visual criterion for the decomposition of a grove]}
We denote by $\Ng^{\infty}:= \{  \emptyset \} \cup \cup_{n \geq 0} \Ng^n$ and by $\Ng^{\infty}_*:=\cup_{n > 0} \Ng^n$.
We present a formal algorithm which recomposes a given grove in terms of binary trees.\\
\textbf{Input:} A grove denoted by $\vec{G} \in \Ng^{\infty}$.\\
\textbf{Output:} A collection of binary trees denoted by $(\vec{v}^{(j)}_{(i)})_{i \in I_j}, \ j \in J$, where $J$ and the $I_j$ are sets such that:
\begin{equation}
\label{genegro}
\vec{G}:= \cup_{j \in J} \ \dotplus_{i \in I_j} \vec{v}^{(j)}_{(i)} \cup  \vec{G}_0,
\end{equation}
where $\vec{G}_0$ is a grove (which is the most general form of a grove).\\
Use the following formula established in Proposition \ref{dend}. If $(\vec{v}_{(i)})_{1 \leq i \leq n}$ is a collection of names of binary trees, then:
\begin{equation}
\label{genegro1}
\dotplus_{1 \leq i \leq n} \ \vec{v}_{(i)} := \vec{v}_{(1)} \dotplus \vec{v}_{(2)} \dotplus \ldots \dotplus \vec{v}_{(n)}=\cup_{\vec{v}_{(1)} \nearrow \vec{v}_{(2)} \nearrow \ldots \nearrow \vec{v}_{(n)} \leq \vec{u} \leq \vec{v}_{(1)} \nwarrow \vec{v}_{(2)} \nwarrow \ldots \nwarrow \vec{v}_{(n)}} \ \vec{u}.
\end{equation}
If a given grove $\vec{g}$ is just the dendriform addition of $n$ vectors
$(\vec{v}_{(i)})_{1, \ldots n}$, then $\min:=\vec{v}_{(1)} \nearrow \vec{v}_{(2)} \nearrow \ldots \nearrow \vec{v}_{(n)}:=(\vec{v}_{(1)}, v_{(1)} \triangleright \vec{v}_{(2)}, \ldots)$ and $\max:=\vec{v}_{(1)} \nwarrow \vec{v}_{(2)} \nwarrow \ldots \nwarrow \vec{v}_{(n)}:=(\vec{v}_{(1)}, v_{(1)} +  \vec{v}_{(2)}, \ldots)$ have to belong to the grove. The first term has the maximum of 1 (first term in the dendriform sum) and the second one the minimum of 1 (last term in the dendriform sum).
Comparing them give automatically the decomposition of the grove $\vec{g}$ into $n$ vectors. Indeed, observe that  whenever the coordinate $\min_i:=1$, $\max_i \in \{1, 1+ v_1, 1+ v_1+v_2, \ldots\}$. By localising the jumps, we determine easily the the lengths of the $\vec{v}_{(i)}$ and thus the $\vec{v}_{(i)}$. Once the $\vec{v}_{(i)}$ obtained, recompute the sum and compare to the grove $\vec{g}$.\\
In general, the grove may be a union of several dendriform sums.  Applying Formula (\ref{genegro1}), observe that every vector of a given dendriform sum starts with $(\vec{v_1}, \ldots)$.\\ 
\textbf{Step 1:} Fix a grove $\vec{G}$ and gather vectors starting with the same 
name of a binary tree.\\
\textbf{Step 2:} Once the kernels are done, take one starting with say $(\vec{v_1}, \ldots)$. Discard $\vec{v_1}$. Subtract the quantity $v_1$ to any element of this kernel and replace negative numbers by 1.\\ 
\textbf{Step 3:} Gather vectors starting with the same 
name of a binary tree, say $(\vec{v_2}, \ldots)$. Reapply Step 2 up to recover all the names of binary trees part of the different dendriform sums.\\
\textbf{Step 4:} Once the different vectors composing the dendriform sums are obtained. Compute these sums and compare with the given grove $\vec{G}$.\\
\textbf{Step 5:} Proceeding that way, it can occur that a dendriform sum give other vectors that those composing $\vec{G}$. That kernel of vectors cannot be replaced by a dendriform sum and is rejected in the notation $\vec{G}_0$. It can also occur that a dendriform sum give part of vectors composing $\vec{G}$. In this case, that part of the grove $\vec{G}$ can be explicitely written in terms of a sum.
\begin{exam}{}
Consider the grove $\vec{G}:=(1,1) \cup (1,2)$. These vectors start with $\vec{v}_1:=(1)$. We discard it and obtain $(1)$ and $(2)$. Substracting 
$v_1:=1$ and replacing negative numbers by 1, leads  the following names of the trees $(1)$ and $(1)$. They both start with $(1)$ which gives us  $\vec{v}_2:=(1)$. We compute the dendriform sum $(1) \dotplus (1)$ and compare it to the given grove we started with and find:  $(1,1) \cup (1,2)= (1) \dotplus (1)$.
\end{exam} 
\begin{prop}
Fix a grove $\vec{G}$. Then, there exists a unique
 collection of binary trees denoted by $(\vec{v}^{(j)}_{(i)})_{i \in I_j}, \ j \in J$, where $J$ and the $I_j$ are sets, (maybe empty)  and a unique grove $\vec{G}_0$, such that:
\begin{equation*}
\label{genegro}
\vec{G}:= \cup_{j \in J} \ \dotplus_{i \in I_j} \vec{v}^{(j)}_{(i)} \cup  \vec{G}_0,
\end{equation*}
\end{prop}
\Proof
Apply the previous algorithm and observe the uniqueness of the vectors obtained by this algorithm.
\eproof
\NB \textbf{[Solving equations]}
Proceeding that way, solutions of first degree equations with unknown can be solved, like for instance,
$ \vec{v} \dotplus \vec{X} := \cup_i \vec{v}_i,$
where only  $\vec{X}$ is not known.

\noindent
As expected, the dendriform addition splits into two binary operations on groves $\dashv$ (Left operation) and  $\vdash$ (Right operation) (pictorially the sign $\dotplus$ gives the signs $\dashv$ and  $\vdash$).
For all $\vec{v} \in \N^n$ and $\vec{w} \in \N^m$,
\begin{eqnarray*}
\vec{v} \dashv \vec{w} &:=& \vec{v}_l \vee (\vec{v}_r \dotplus \vec{w}):= \bigcup_{\vec{v}_l \vee (\vec{v}_r \nearrow \vec{w}) \leq \vec{u} \leq \vec{v} \nwarrow \vec{w}} \ \vec{u}\ \ \textrm{when} \ \ \vec{v} \not=(0), \\
\vec{v} \vdash \vec{w} &:=& (\vec{v} \dotplus \vec{w}_l) \vee \vec{w}_r:= \bigcup_{\vec{v} \nearrow \vec{w} \leq \vec{u} \leq (\vec{v} \nwarrow \vec{w}_l) \vee \vec{w}_r  } \ \vec{u}, \ \ \textrm{when} \ \ \vec{w} \not=(0), \\
\end{eqnarray*}
These operations are extended to $\Ng^\infty$, that is to groves and to $(0)$ by distributivity with respect to the disjoint union and verify the axioms:
$(\vec{u} \dashv \vec{v} ) \dashv \vec{w} = \vec{u} \dashv (\vec{v} \dotplus  \vec{w}), \ 
(\vec{u} \vdash \vec{v} )\dashv w = \vec{u} \vdash(\vec{v} \dashv \vec{w}), \ 
(\vec{u} \dotplus \vec{v} )\vdash \vec{w} = \vec{u} \vdash(\vec{v} \vdash \vec{w}), \  $
$(0) \dashv \vec{v}:= \emptyset =: \vec{v} \vdash (0)$ 
and $(0) \vdash \vec{v}= \vec{v}=\vec{v} \dashv (0)$. We set
$\emptyset \circ \vec{v}:=\vec{v} \circ \emptyset := \emptyset$, for 
$\circ \in \{\dashv, \ \vdash\}$ and any grove $\vec{v}$. The r\^ole of the empty set will be explained later.
The symbols $(0) \dashv (0)$ and $ (0) \vdash (0)$ are not defined, though $(0) \dotplus (0):=(0).$ Moreover, $(\vec{u} \vdash \vec{v})^\dagger:= \vec{v}^\dagger \dashv \vec{u}^\dagger$ and $(\vec{u} \dashv \vec{v})^\dagger:= \vec{v}^\dagger \vdash \vec{u}^\dagger$.
The tricks for computations are the following.
$ (\vec{u},1) \dashv \vec{v}:= \vec{u} \vee \vec{v}:=(\vec{u},1, u+1+ \vec{v}),$
and
$\vec{u} \vdash (1,\vec{v}):=\vec{u} \vee \vec{v}.$  
\subsubsection{The dendriform involution}
\label{dendinvosubsect}
There is an involution on $\Ng^\infty:=\cup_{n \geq 0} \Ng^n$ denoted by $\dagger$ and defined by
$ (\vec{v} \vee \vec{w})^\dagger:= \vec{w}^\dagger \vee  \vec{v}^\dagger.$ That is $ (\vec{v}, 1, v+\vec{w})^\dagger:=(\vec{w}^\dagger, 1, w+\vec{v}^\dagger)$.
Doing so, observe that $ (\vec{v} \nearrow \vec{w})^\dagger=\vec{w}^\dagger \nwarrow \vec{v}^\dagger$ and
$ (\vec{v} \nwarrow \vec{w})^\dagger=\vec{w}^\dagger \nearrow \vec{v}^\dagger$. Therefore,
$ (\vec{v} \dotplus \vec{w})^\dagger= \vec{w}^\dagger \dotplus \vec{v}^\dagger,$
\textit{i.e.}, $(\Ng^\infty, \dotplus, \dagger)$ is an involutive graded monoid. Observe that $(1)^\dagger:=(1)$ and by convention, we set $(0)^\dagger:=(0)$. We now state some properties of the involution on trees by giving another representation of the Catalan numbers.
\begin{prop}
Fix $n \geq 1$ and let $Inv[n]:=\{\vec{v}  \in \N^n, \ \vec{v}:=\vec{v}^\dagger \}$. Then, $Inv[2n]= \emptyset$ and $\card(Inv[2n+1])=c_n:=\frac{1}{n+1}(^{2n}_{\ n})$.
\end{prop}
\Proof
Observe that $\vec{v}:=\vec{v}^\dagger$ if and only if there exists a unique $\vec{w}$ such that $\vec{v}=\vec{w} \vee \vec{w}^\dagger$.
\eproof

\noindent
\textbf{[Trick to name $\vec{v}^\dagger$.]}
Fix $\vec{v} \in \N^n$. There exists a very simple way to name
$\vec{v}^\dagger$. Associate with $\vec{v}$, its complete expression in $\bra x_1, \ldots, x_n,x_{n+1}, (,) \ket$. Relabel $x_{n+1}$ by $x_{1}$,
$ x_n$ by $x_{2}$ and so on. Read therefore from left to right such a monomial. The vector $\vec{v}^\dagger \in \N^n$ is obtained from the following construction. The coordinate $v^\dagger_i:=i$, for all $1 \leq i \leq n$, if and only if there is a $)$ at the right hand side of $x_i$ and $v^\dagger_i:=j$ if the most left parenthesis $($ at the left hand side of $x_i$ closes a $)$ open in $x_j$. This works since the involution
on binary trees is a symmetry with regards to the root axis, which can also be viewed as a symmetry with regards to an axis perpendicular to it --the Mirror axis-- giving then the mirror image of the tree and thus its involution.
\begin{center}
\includegraphics*[width=6cm]{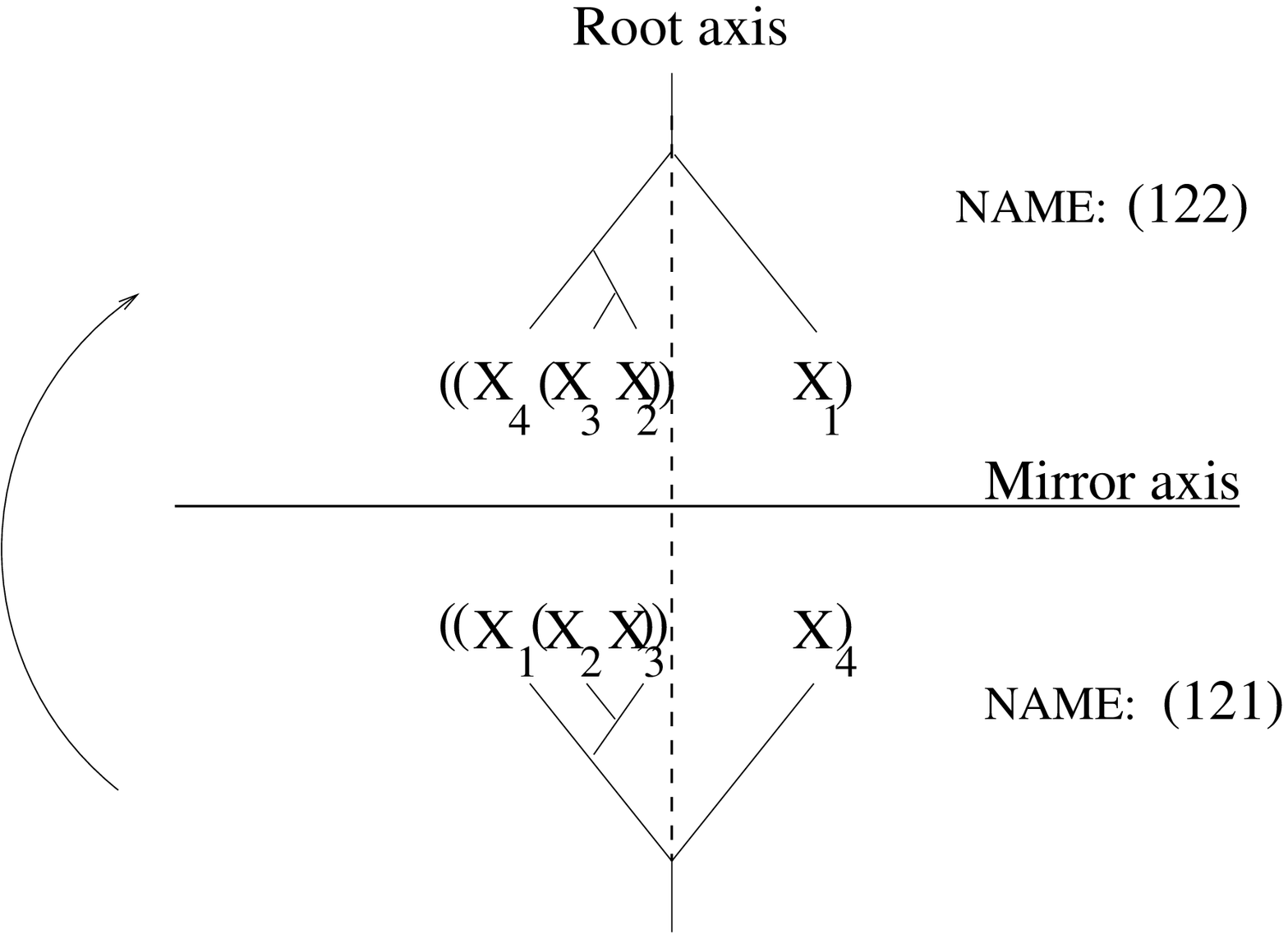}
\end{center}
\begin{prop}[Lattice anti-automorphism.]
\label{invinf}
Let $\vec{v}, \vec{w} \in \N^n$. Then, the dendriform involution is a lattice anti-automorphism, \textit{i.e.},
$ \vec{v} <  \vec{w} \Leftrightarrow \vec{w}^\dagger <  \vec{v}^\dagger.$
Consequently, $M(\vec{v},\vec{w})=M(\vec{w}^\dagger,\vec{v}^\dagger)$, for any names
of trees.
\end{prop}
\Proof
Fix $\vec{v}, \vec{w} \in \N^n$ with $\vec{v} <  \vec{w}$. 
We will check the case when both  $v_i=i=w_i$. In this case $v^\dagger_{n+1-i} =j$ and $w^\dagger_{n+1-i} =j'$ with $j' \leq j < N+1-i.$ Indeed, suppose the existence of a $)$, most external parenthesis standing at the right hand side of $x_{j'}$ and closing one $($
open in $x_i$ in the complete expression associated with $\vec{w}$. As $\vec{v} <  \vec{w}$, we get $v_{j'} \leq w_{j'}:=i$. If $v_{j'}:=k<i$, then this means that the most external parenthesis $)$ standing at the right hand side of $x_{j'}$ in the complete expression associated with $\vec{v}$ closes one $($ open in $x_k$. This implies that $v^\dagger_{n+1-i} \leq w^\dagger_{n+1-i}$. Checking every possibility leads to the conclusion that $v^\dagger_{i} \leq w^\dagger_{i}$ for all $1 \leq i \leq n$. The proof is complete since the dendriform involution is an involution.
For the last claim,
the dual lattice of $(\N^n,<)$ with order $\leq^*$ is such that
$\vec{v}  \leq^* \vec{w} \Leftrightarrow \vec{w}  \leq \vec{v} $. Therefore,
$\vec{v}  \leq^* \vec{w} \Leftrightarrow \vec{v}^{\dagger}  \leq \vec{w}^{\dagger}$, for any vectors of $\N^n$. The last claim holds since
$M^*(\vec{v}, \vec{w}):=M(\vec{w}, \vec{v})$ (see Prop. 3 p345 of \cite{Rota}).
\eproof
\subsubsection{The dendriform multiplication}
The following idea developed by J.-L. Loday consists to replace the polynomial ring $K[X]$,  (basis $(X^n)_{n \in \Na}$) and well-known equations $X^nX^m:=X^{n+m}$ and $(X^n)^m:=X^{nm}$ related to the usual arithmetic on $\Na$,
by planar binary trees. Instead of writting $K[X]$, one could have chosen $K[\Na]$ to denote this polynomial ring.
Consider the $K$-vector space $K[\Ng^\infty_*]$ spanned by the basis $\{X^{\vec{v}}, \ \vec{v} \in \N^\infty_*\}$. The space $K[\Ng^\infty_*]$ has a natural
dendriform algebraic structure given by:
$ X^{\vec{u}} \prec X^{\vec{v}} := X^{\vec{u} \dashv \vec{v}}$ and  $X^{\vec{u}} \succ X^{\vec{v}} := X^{\vec{u} \vdash \vec{v}},$
with the convention: $X^{\vec{u} \cup \vec{v}}:=X^{\vec{u}} + X^{\vec{v}}.$ As expected,
$ X^{\vec{u}} \star X^{\vec{v}}:=X^{\vec{u} \dotplus \vec{v}},$ 
where $\star$ is the associative product, sum of $\prec$ and $\succ$. This nonunital associative algebra, another representation of the free dendriform dialgebra on one generator, here $X^{(1)}$, can be augmented by adding the unit $1:=X^{(0)}$ so that, $K[ \Ng^{\infty}]:=K[\Ng^{\infty}_*] \oplus K \cdot 1.$ By convention, we set $X^{\emptyset}=0$.
As usual, the operations $\prec$ and $\succ$ can be partially extended to $K[ \Ng^{\infty}]$ by declaring that: $1 \succ X^{\vec{v}} :=X^{\vec{v}}=: X^{\vec{v}} \prec 1$, for $\vec{v} \not=(0)$ and vanish otherwise, explaining the presence of the empty set. For instance, $1 \prec X^{\vec{v}} := X^{(0) \dashv \vec{v}}:= X^{\emptyset} := 0$, as expected.
The notation $K[\Na]:=K[X]$ stands for the usual polynomial algebra on one variable say $X$. As $X^nX^m:= X^{n+m}$ and $\Na$ is invariant by addition, one can use also the notation $K[\Na]$ without any ambiguity. However, $K[\N^{\infty}]$ is not invariant by the dendriform addition, that is why we choose the notation $K[ \Ng^{\infty}]$ and not $K[\N^{\infty}]$. The notation $K[\Ng^{\infty}]$ stands for the $K$-vector space spanned by $\{X^{\vec{v}}, \  \vec{v} \in \N^{\infty}\}$.
\begin{defi}{[Dendriform multiplication \cite{Lodayarithm}]}
The \textit{dendriform multiplication} $\ltimes: \Ng^n \times \Ng^m \xrightarrow{} \Ng^{nm}  $  is given by
$\vec{u} \ltimes \vec{v}:= \omega_{\vec{u}}(\vec{v}),$
for all $\vec{u}$ and $\vec{v}$, names of binary trees and extended to groves \textit{via} distributivity on the left with respect to the disjoint union, \textit{i.e.}, $(\vec{u} \cup \vec{v})\ltimes \vec{w}:=\vec{u}\ltimes \vec{w}  \cup \vec{v} \ltimes \vec{w}.$
\end{defi}
For instance, as $(1,2):=(1) \dashv (1)$, we get $(1,2) \ltimes \vec{v}:=(\vec{v}) \dashv (\vec{v})$. Therefore, $(1,2) \ltimes (1,1):=(1,1) \dashv (1,1):=(1,1,3,3)$.
The dendriform multiplication is associative, not commutative, distributive on the left with regards to the dendriform addition $\dotplus$, has the neutral element $(1)$ and is compatible with the involution $\dagger$, $(\vec{u}  \ltimes \vec{v})^\dagger=\vec{u}^\dagger  \ltimes \vec{v}^\dagger$. Moreover, the neutral element for $\dotplus$, \textit{i.e.}, $(0)$, is by convention a left anihilator for $\ltimes$, \textit{i.e.,}
$(0)  \ltimes \vec{u}:=(0)$.
A vector $\vec{w} \in \N^{n}$ is said to be \textit{prime} if there exists no vector $\vec{v} \in \N^{m}$ and $\vec{v'} \in \N^{m'}$, with $n=mm'$ such that $ \vec{w}:=\vec{v} \ltimes \vec{v'}$. In general, the dendriform product of two vectors gives a grove. However, observe there are two unique ways  to obtain a vector.  The first one is to consider $(1) \dashv (1) \ltimes \vec{v}$, with $\vec{v}:= \vec{v}_1 \vee (0) $ and the second one is to consider $(1) \vdash (1) \ltimes \vec{v}$, with $\vec{v}:= (0) \vee \vec{v}_1$. In the first case, we obtain $(1) \dashv (1) \ltimes \vec{v}:= \vec{v}_1 \vee (\vec{v}_1 \vee (0) )$ and in the second case,
$(1) \vdash (1) \ltimes \vec{v}:= ((0) \vee \vec{v}_1) \vee \vec{v}_1.$ We summarize our discussion by the following proposition.
\begin{prop}
Any vector of $\N^{2n+1}$ is prime for the arithmetree just described.
Whereas, there exist $2c_{n}$ nonprime vectors in $\N^{2n+2}$. They are of the forms:
$ ((0) \vee \vec{v}) \vee \vec{v} \ \ \textrm{and} \ \  \vec{v} \vee (\vec{v} \vee (0)),$
with $\vec{v} \in \N^{n}$.
\end{prop}
\begin{prop}{[Right and left cancellations]}
Let $\vec{v} \in \N^{n}, \vec{u},\vec{w} \in \N^{m}$. Then,
$$ \vec{v} \ltimes \vec{u} = \vec{v} \ltimes \vec{w} \Leftrightarrow \vec{u}=\vec{w}, \ \
 \vec{u} \ltimes \vec{u} = \vec{w} \ltimes \vec{w} \Leftrightarrow \vec{u}=\vec{w}, \ \
 \vec{u} \ltimes \vec{v} = \vec{w} \ltimes \vec{v} \Leftrightarrow \vec{u}=\vec{w}.$$
\end{prop}
\Proof
The first claim is obtained by observing that the first operation appearing in $\omega_ {\vec{v}}((1))$ is either $\vdash$ or $\dashv$. Therefore, in both cases, the vectors composing the groves $\vec{v} \ltimes \vec{u}$ and $\vec{v} \ltimes \vec{w}$ will start with $(\vec{u}, \ldots)$, resp. with $(\vec{w}, \ldots)$. The same remark applies also for the second claims. To complete the proof, observe that the dendriform multiplication acting on the right hand side is the unique dendriform automorphism which maps the generator $X^{(1)}$ to $X^{\vec{v}}$ in $K[\Ng]^{\infty}$.
\eproof

\noindent
Recall that the free dendrifrom algebra on the generator $X^{(0)}$ is
linked to the free associative $L$-algebra on the same generator, the operations being given by the under and over operations $\nearrow$ and $\nwarrow$. 
\begin{prop}
Let $\vec{u}$ be a name of a binary tree.
Then, $\varpi_{\vec{u}}((1))$ can be obtained from $\omega_{\vec{u}}((1))$ by replacing the symbols $\vdash$ by $ +_{_{\nearrow}}$ and $\dashv$ by $+_{_{ \nwarrow}}$. We name the middle term the vector so obtained.
If a grove is not prime for the dendriform arithmetics, then its middle term will be not prime for the $L$-arithmetics.
\end{prop}
\Proof
Proceed by induction. It is true for $n:=1,2,3$ (checked by hand).
Observe that $\vec{u}_l \vdash (1) \dashv \vec{u}_r= \vec{u}_l \vee \vec{u}_r= \vec{u}_l +_{_{\nearrow}} (1) +_{_{\nwarrow}} \vec{u}_r$. Therefore,
$\omega_{\vec{u}_l}(1) \vdash (1) \dashv \omega_{\vec{u}_r} (1)$ gives
$\varpi_{\vec{u}_l}(1) +_{_{\nearrow}}  (1) +_{_{\nwarrow}} \varpi_{\vec{u}_r} (1)$ by replacing the symbols $\vdash$ by $+_{_{\nearrow}}$ and $\dashv$ by $+_{_{ \nwarrow}}$.
\section{Bijection between noncrossing partitions and binary trees}
We recall a bijection between noncrossing partitions and binary trees. A noncrossing partition of the set $\{ 1,2,3, \ldots,n, \}$ is a decomposition $\pi:=\{V_1, \ldots ,V_r \}$ of $S$ into disjoint and nonempty sets $V_i$, called blocks, such that for all $1 \leq p_1, q_1,  p_2, q_2 \leq n,$ the following does not occur:
there exist
$ 1 \leq p_1 < q_1 < p_2 < q_2$ with $p_1 \backsim_{\pi} p_2  \nsim_{\pi} q_1 \backsim_{\pi} q_2,$
where for all $1 \leq p,q \leq n,$ $p \backsim_{\pi} q$ means
that $p$ and $q$ belong to the same block of $\pi$.
The set of noncrossing partitions made out of the elements $1,2,3, \ldots,n$ is denoted by $NC(n)$ .  In low dimensions, these sets are,
\begin{center}
\includegraphics*[width=13cm]{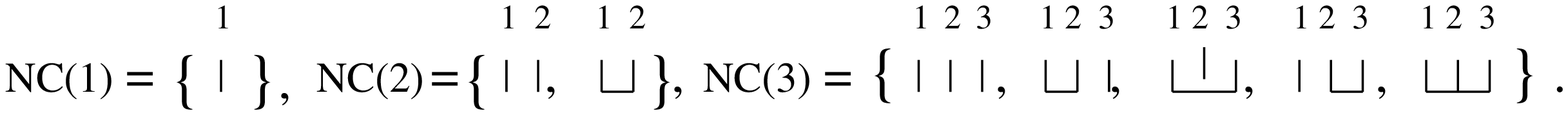}
\end{center}
There is a natural poset structure given by the refinement order.
In the sequel, an interval of a bloc $V$ is a sequence of numbers all linked one another. Every bloc can be decomposed uniquely in several intervals.
A bijection between noncrossing partitions and binary trees is determined by the following algorithm in 2 steps.
\begin{enumerate}
\item {Let $\tau \in Y_n$ be a $n$-tree, $n>0$. As the tree is planar and binary, the notion of left and right has still a meaning.  As the tree is rooted, denote by 1 the root. This gives a Cartesian plan of dimension two denoted by $(1, R, L)$ where the axis $R$ (resp. $L$) is the line passing through 1 and identified with the most right (resp. left) branch. Pictorially, we get:
\begin{center}
\includegraphics*[width=2cm]{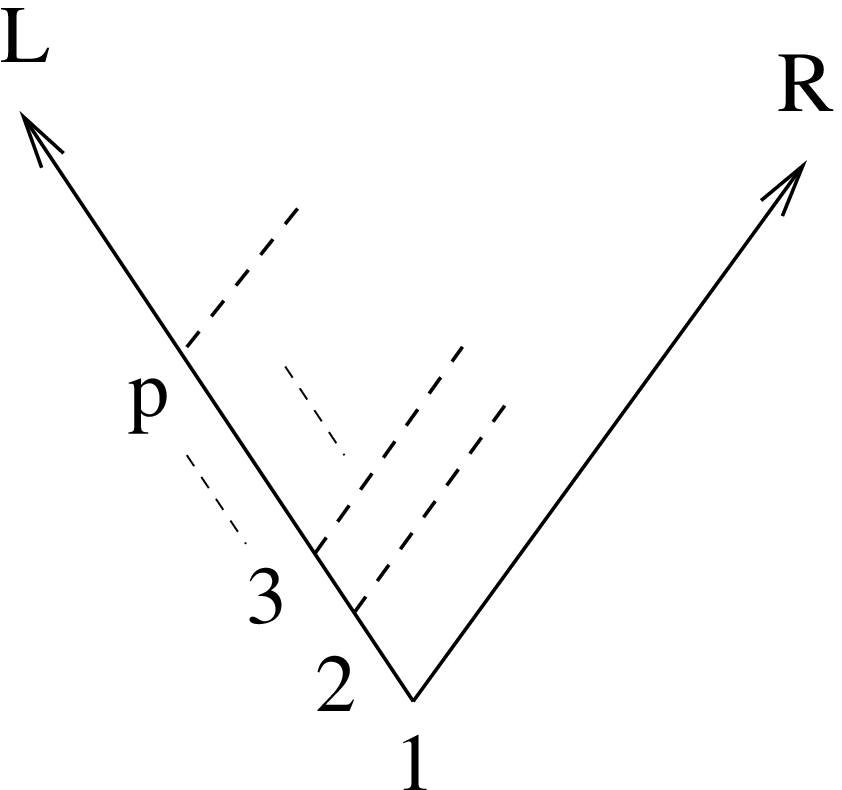}
\end{center} 
Starting with the origin of the Cartesian plan here $(1,L,R)$, \textit{i.e.}, with 1, increment of a unit all the $p>1$ branchs linked to the axis $L$, this gives $1,2,3 \ldots,p$. If there is no branch at the left hand side of 1, give 2 to the closest vertex at the right hand side of 1 and reapply the algorithm.}
\item {Once arrived in the vertex $p$. If there is a vertex to the right of $p$, give the number $p+1$ to it and reapply 
the algorithm in the Cartesian plan $(p+1, L, R)$ modelling now the subtree with root the vertex $p+1$. If not, go to the vertex $p-1$ and reapply the algorithm at Step 2.}
\end{enumerate} 
Once all vertices of the tree are labelled, a unique noncrossing partition is obtained by the following trick. Put a vertical segment under each numero $1,2, \ldots , n$ and
link $p$ to $q>p$ if  $q$ is the closest vertex at the right hand side of $p$. One can view this partition 
as the `projection' ---(by abuse of language)--- parallel to the axis $L$ of all the branches of the trees on the axis $R$ in the Cartesian plan $(1,L,R)$ if the branches are all drawn either parallel to the axis $L$ or $R$.
Here is an example.
\begin{center}
\includegraphics*[width=6cm]{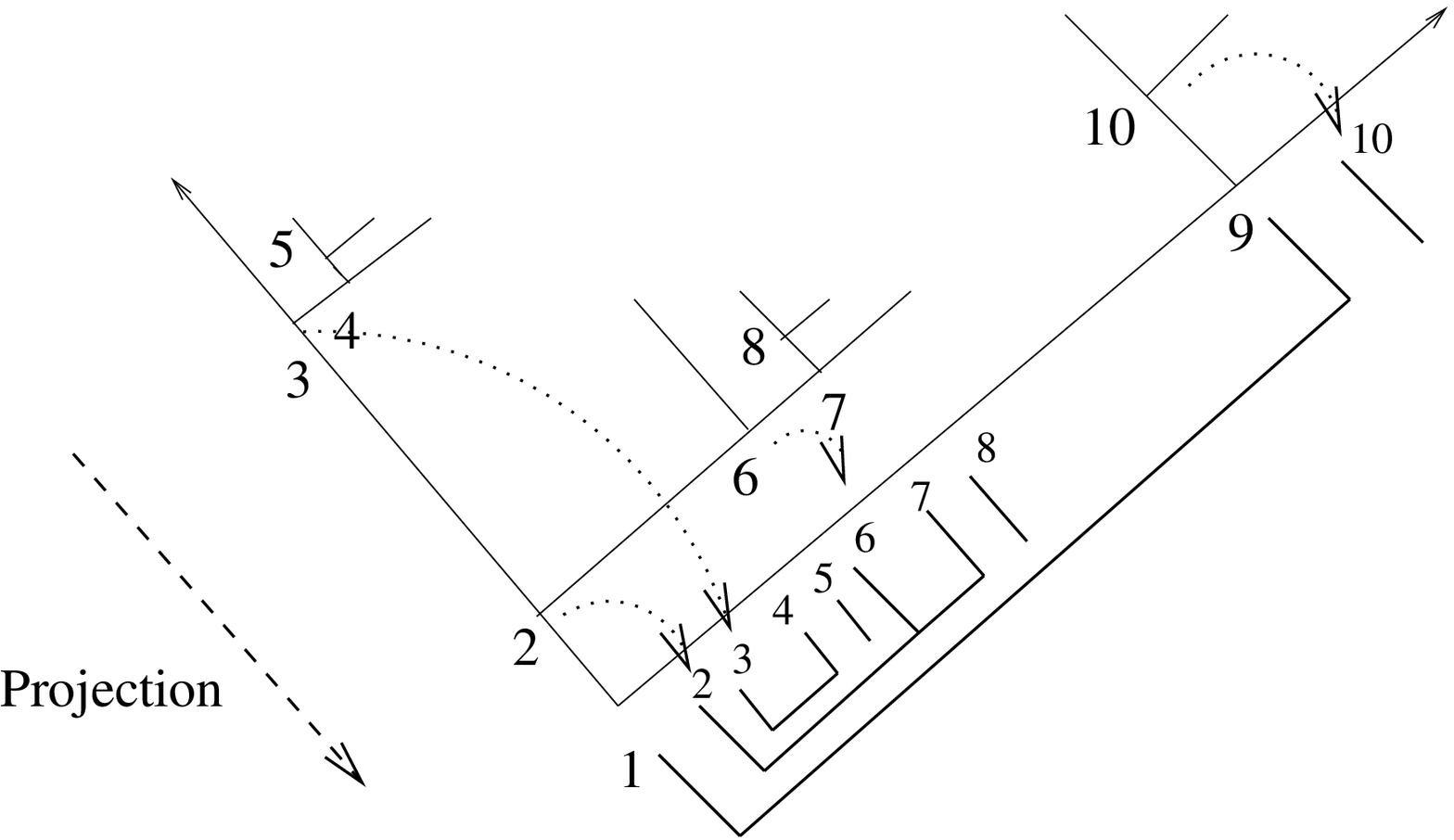}
\end{center} 
To recover the binary tree from its noncrossing partition, proceed as follows. For pedagogical reason, we will proceed on the example just above. By construction $1$ has to be the root of the binary tree. Therefore, draw the Cartesian plan $(1,R,L)$. The root is linked to $9$ so $(9,R,L)$ has to denote the closest vertex (with 2 branchs) at the right hand side of 1. However 1+1=2 is not linked to 1, so $(2,R,L)$ is at its left hand side. Now focus only on the numbers between 2 and 8. This gives another noncrossing partition. Reapply then the previous algorithm by asking who is on the right and on the left hand side of a given number. Observe for instance that 5 is a leaf since there is no element at its right hand side (5 is not linked to any number) and there is no element at its left hand side neither since 5+1=6 is linked to 2 and so cannot be at the left hand side of 5. 

\NB
The previous algorithm gives a bijection $\Pro: Y_n \xrightarrow{} NC(n)$, for all $n>0$. The set $NC(n)$ can be equipped with the Tamari order and $Y_n$ with the refinement one.
There exists another way to associate bijection to planar binary trees (compare to \cite {Loday}). Indeed, any noncrossing partition models naturally a bijection written in disjoint cycles. The bijection modelled by our previous example is  $(1,9)(10)(2,6,7)(8)(3,4)(5)$. 
There exists an involution on $NC(n)$ induced by the one introduced on $Y_n$. If $\pi$ denotes such a partition, then $\pi^\dagger$ can be easily constructed by the projection on the axis $L$ in the Cartesian plan $(1,R,L)$ of the mirror image of the tree associated with $\pi$.
This bijection will play an important r\^ole for the sequel of this paper. It will give \textsf{NCP}-operads and a connection with the free dendriform dialgebra on one generator.
\section{A reformulation of free probability}
There are several kind of geometry available, among them, the most well-known being the Euclidean geometry. The same thing holds in Probability theory where the most well-known is of course the classical probability, \textit{i.e.}, defined by the Kolmogorov axioms leading thus to the usual stochastic independence. However, the introduction of quantum mechanics in physics has paved the way to other challenging stochastic independences. Among these ones, lies a noncommutative probability theory equipped with the so-called free stochastic independence. It is rooted in $C^*$-algebras and has been pointed out first by D. Voiculescu. Later, a complementary point of view was given by R. Speicher \cite{Speicher}, inspired by previous works of G.-C. Rota \cite{Rota2}. 
\subsection{Action of arithmetree on $\B-\B$-bimodule and operads}
In the sequel, $\B$ denotes a unital associative algebra (most of the time a unital $C^*$-algebra for applications) and $\M$ is a $\B-\B$-bimodule. We denote by $\M^{\otimes_{\B} n}$, the space $\M \otm \M \otm \ldots \otm \M$, $n$ times and by convention $\M^{\otm 0}:=B$.
By abuse of language and sometimes to ease notation, we will use equivalently trees and/or their names.
One of the aims of this  part is to describe the action of the space $K[\YY_{\infty}]$  --- or equivalently $K[\Ng^{\infty}]$ --- defined in Section 3 and equipped with its arithmetree onto the bimodule $\M$. 
\subsubsection{$\NCP(\B)$-Operads}
To introduce the action of binary trees \textit{via} in terms of noncrossing partitions, we will need the concept of 
noncrossing partitions operads,  $\NCP(\B)$-operads for short. This concept is inspired from \cite{Speicher} and is different, though similar, from
the regular $K$-linear operads definition see \cite{Lodayren} or the introduction of this paper.
\begin{defi}{}
Let $\B$ be an associative $K$-algebra.
A \textit{$\NCP(\B)$-operad} $\P$ (without unit) over a $\B-\B$-bimodule $\M$ is the data of a family of finite dimension $K$-vector spaces $(\P(n))_{n >0}$, whose basis elements $\mu$ are $\B-\B$-bimodule $n$-ary operations with values in $\B$, \textit{i.e.},
$\mu: \M^{\otm n} \xrightarrow{} \B,$
and equipped with a family of composition maps ($(\circ_i)_{i>0}$) verifying the following relations,
\begin{enumerate}
\item {For all $\mu \in \P(m)$ and $\nu \in \P(n)$ and $1\leq i \leq m+1$,  $\mu \circ_i \nu \in \P(m+n)$.}
\item {For all $\lambda \in  \P(l)$, $\mu \in \P(m)$ and $\nu \in \P(n)$, 
$$ (\lambda \circ_i \mu)\circ_{j+m} \nu=(\lambda \circ_j \nu)\circ_i \mu, \ \ 1\leq i \leq j \leq l+1,$$
$$ \lambda \circ_i (\mu \circ_{j} \nu)=(\lambda \circ_i \mu)\circ_{i+j-1} \nu, \ \ 1 \leq i \leq l+1, \ \ \ 1 \leq j \leq m+1.$$}
\end{enumerate} 
An \textit{augmented} $\NCP(\B)$-operad $\P^{+}$ is the data of a $\NCP(\B)$ -operad $\P$ such that $\P^{+}(n):=\P(n)$, for $n>1$ and 
$\P^{+}(1):=K \oplus \P(1)$.
\end{defi}
A noncrossing partition $\mu \in NC(n)$ is said to be decorated by a set $Col$ if a unique color of $Col$ is associated
with each interval composing it. 
Observe that decorated noncrossing partitions give special decorated binary trees, \textit{i.e.}, binary trees whose all vertices of a SW-NE-branch have the same color. 
We now give an example of such $\NCP(\B)$-operad by mixing results in Section 3 on the free dendriform dialgebra
on one generator and ideas developed from noncommutative probability. 
\subsubsection{$\NCP(\B)$-Operads and dendriform stucture on $\B-\B$-bimodules}
Let $\M$ be a $\B-\B$-bimodule. The dendriform dialgebra over $K[\Ng^{\infty}_*]$ induces a dendriform dialgebra structure on the following $K$-vector space:
$ \textsf{Dend}_{\B}(\M):=\bigoplus_{n \geq 1} \ K[\Ng^n] \otimes \M^{\otm n},$
by declaring that:
\begin{eqnarray*}
X^{\vec{v}} \otimes \kappa \prec X^{\vec{w}} \otimes \kappa' &:=& X^{\vec{v}}\prec X^{\vec{w}} \otimes \kappa\kappa', := X^{\vec{v} \dashv \vec{w}} \otimes \kappa\kappa',\\
X^{\vec{v}} \otimes \kappa \succ X^{\vec{w}} \otimes \kappa' &:=& X^{\vec{v}}\succ X^{\vec{w}} \otimes \kappa\kappa':= X^{\vec{v} \vdash \vec{w}} \otimes \kappa\kappa',
\end{eqnarray*}
for any tensor $\kappa \in  \M^{\otm n}$ and $\kappa' \in  \M^{\otm n'}$.
Of course, to make operations well-defined, we have to divide out by the following relations, $X^{\vec{v}} \otimes \kappa b \prec X^{\vec{w}} \otimes \kappa'=X^{\vec{v}} \otimes \kappa  \prec X^{\vec{w}} \otimes b \kappa'$
and $X^{\vec{v}} \otimes \kappa b \succ X^{\vec{w}} \otimes \kappa'= X^{\vec{v}} \otimes \kappa  \succ X^{\vec{w}} \otimes b\kappa'$, for all $b \in B$.
We still refer to $ \textsf{Dend}_{\B}(\M)$ to be the free dendriform dialgebra over $\M$. 
Inspired by \cite{Speicher}, we define a family of $n$-ary operations, 
$ f^{(n)}: \M^{\otm n} \xrightarrow{} \B,$
which are $\B-\B$-bimodule maps, \textit{i.e.}, $f^{(n)}(ba_1 \otm a_2 \otm \ldots \otm a_nb'):=bf^{(n)}(a_1 \otm a_2 \otm \ldots \otm a_n)b'$, for all $n>0$ and $b,b' \in \B$ and $a_1, \ldots , a_n \in \M$.
With each family $(f^{(n)})_{n \geq 1}$, we associated the following operator valued function:
$$ \hat{f}:= (f^{(n)})_{n \geq 1}: \DM \xrightarrow{} B, \ \ \
(X^{\vec{v}} \otimes a_1 \otm a_2 \otm \ldots \otm a_n) \mapsto \hat{f}(X^{\vec{v}} \otimes a_1 \otm a_2 \otm \ldots \otm a_n),$$
defined \textit{via} the following recursive prescrition:
With any monomial $X^{\vec{v}}$, is associated a unique noncrossing partition $\Pro(\vec{v})$, constructed from the algorithm described in the previous section. 
Identify this partition to the tensor $a_1 \otm a_2 \otm \ldots \otm a_n$. 
Localise the most nested block of lenght $p \leq n$ and apply
the $p-ary$ operations, giving thus an operator in $\B$. Then,
reapply this procedure. In the sequel, we will write:
$$ (X^{\vec{v}} \otimes a_1 \otm a_2 \otm \ldots \otm a_n) \mapsto \hat{f}(X^{\vec{v}} \otimes a_1 \otm a_2 \otm \ldots \otm a_n):= \hat{f}(\Pro(\vec{v}) \rightsquigarrow (a_1 \otm a_2 \otm \ldots \otm a_n)),$$
to denote that action of the noncrossing partition $Pr(\vec{v})$.
The following examples will be better than a fastidious description. Here are three examples (recall that $a \otm ba':=a b \otm a'$, for $b \in \B$):
\begin{enumerate}
\item {$\hat{f}(X^{^{^{\treeACA}}} \otimes a_1 \otm a_2 \otm a_3):= f^{(2)}(a_1 \otm f^{(1)}(a_2) \otm a_3):=f^{(2)}(a_1 f^{(1)}(a_2) \otm a_3).$}
\item {$\hat{f}(X^{^{^{\treeCAB}}} \otimes a_1 \otm a_2 \otm a_3):= f^{(2)}(a_1 \otm a_2) \otm f^{(1)}(a_3):=f^{(2)}(a_1 \otm a_2)  f^{(1)}(a_3).$}
\item {Let $\Pro(\vec{v}):=(1,9)(2,6,7)(3,4)(5)(8)(10)$ be  the noncrossing partition represented in Section~2 and get:
$$ \hat{f}(X^{\vec{v}} \otimes a_1 \otm a_2 \otm \ldots \otm a_{10}) := \hat{f}(\Pro(\vec{v}) \rightsquigarrow (a_1 \otm a_2 \otm \ldots \otm a_{10})):=$$
$$ f^{(2)}(a_1 \otm f^{(3)}(a_2 \otm f^{(2)}(a_3 \otm a_4) \otm f^{(1)}(a_5) \otm a_6 \otm a_7) \otm f^{(1)}(a_8) \otm a_9) \otm f^{(1)}(a_{10}), $$
and obtain,
$ f^{(2)}(a_1 f^{(3)}(a_2 f^{(2)}(a_3 \otm a_4)  f^{(1)}(a_5) \otm a_6 \otm a_7) \otm f^{(1)}(a_8)  a_9)  f^{(1)}(a_{10}).$}
\end{enumerate} 
\NB
Proceeding that way, observe that $f^{(n)}( a_1 \otm a_2 \otm \ldots \otm a_{n})$ and the action of the maximal element $\textbf{1}_n $ of $NC(n)$ on the $n$-tensor, \textit{i.e.}, $f^{(n)}( \textbf{1}_n  \rightsquigarrow a_1 \otm a_2 \otm \ldots \otm a_{n}),$ coincide.
\NB
We can slightly reformulate this framework using the concept of $\NCP(\B)$-operad. Set $Col:=\{f^{(n)}, \ n>0 \}$ be the color set made out of the $n$-ary operations $f^{(n)}$. Observe that with each noncrossing partition, a unique
decorated noncrossing partition can be associated. Introduce
the object $\P[\hat{f}]$ made out of a family of the $K-$vector spaces $(\P[\hat{f}](n))_{n>0}$ and the family of composition $(\circ_i)_{i>0}$ defined by induction as follows.
The $K$-vector space $\P[\hat{f}](1)$ is spanned by $f^{(1)}$ and $\P[\hat{f}](p:=n+m)$ by
the elements $f^{(p)}$ and $\mu \circ_i \nu$, with 
$\mu \in \P[\hat{f}](m)$ and $\nu \in \P[\hat{f}](n)$ and $1\leq i \leq m+1$, where:
$$\mu \circ_i \nu(a_1 \otm \ldots \otm a_{n+m}):= (a_1 \otm \ldots \otm a_{i-1} \otm \nu (a_i \otm \ldots \otm a_{i+n-1}) \otm a_{i+n} \otm \ldots \otm a_{m+n}),$$
for all $a_1, \ldots, a_n \in \M$ and not in $\B$.
From a noncrossing partition, one can easily write
its action on tensor elements in terms of compositions maps.
The following example will fix ideas.
\begin{exam}{}
Consider again $\Pro(\vec{v}):=(1,9)(2,6,7)(3,4)(5)(8)(10)$,   the noncrossing partition represented in Section~2. 
Read the partition from left to right. Take the first encountered interval, say with $p$ elements ---(here $\{1,9\}$ and $p:=2$)--- and (thus) starting with $1$. Take the second encountered interval, starting with say $n$, and with say $q$ elements ---(here $\{2,6,7 \}$ and $q:=3$)--- and write $f^{(p)} \circ_n f^{(q)} \ldots$. Reapply the algorithm. 
We obtain:
$$\hat{f}(X^{\vec{v}} \otimes a_1 \otm a_2 \otm \ldots \otm a_{10}) :=f^{(2)} \circ_2 f^{(3)} \circ_3  f^{(2)} \circ_5 f^{(1)} \circ_{8} f^{(1)} \circ_{10}  f^{(1)}(a_1 \otm a_2 \otm \ldots \otm a_{10}).$$
\end{exam}  
\begin{theo}
\label{tgb}
Let $\B$ be an associative algebra and $\M$ be a $\B-\B$-bimodule. Let $\hat{f}:=(f^{(n)}:\DM \xrightarrow{} \B)_{n>0}$ be a family of $n$-ary $\B-\B$-bimodule operations.  Then, $(\DM, \hat{f})$ induces a $\NCP(\B)$-operad, $\P[\hat{f}]$, over the $B-B$-bimodule $\M$. 
\end{theo}
\Proof
This theorem summarizes the previous discussions.
\eproof
\NB
The $\NCP(\B)$-operad such obtained can be augmented as follows. By convention, $\Pr(\vec{0})$ is identified with the `noncrossing partition' $\emptyset$. Then, set $\P(1)^{+}:=K \oplus \P(1)$, the $K$-vector space spanned by the set $\{ f^{(1)},  f^{(\emptyset)}  \}$, where $f^{(\emptyset)}(b):=b$, for all $b \in \B$.
\subsection{Cumulants and moments in free probability}
We recall some properties of free probability \cite{Speicher}. Fix $B$, a unital associative algebra. Let $(\M, \phi)$, be a noncommutative probability space, that is a $B-B$-bimodule endowed with a unital associative algebra equipped with a $B-B$-bimodule map $\phi:\M \rightarrow B$ such that $\phi(1)=1$, that is $\phi(b)=b$, for any $b \in B$. Let $\M_1, \ldots, \M_n,$ be $n$ unital $B-B$-subalgebras of $\M$. It is said that $\M_1, \ldots, \M_n,$ are stochastically free if $\phi(a_1\ldots a_n)=0$ under the following conditions. For all $1\leq i \leq n$, $\phi(a_i)=0$ and for $a_1 \in \M_{\epsilon_1}, \ldots, a_n \in \M_{\epsilon_n}$, $\epsilon_1 \not= \epsilon_2$, $\epsilon_2 \not= \epsilon_3$, $\ldots$, $\epsilon_{n-1} \not= \epsilon_{n}$.
It has been shown by R. Speicher, that this definition can be reformulated in terms of noncrossing partitions equipped with the reffinement order. For that, he introduced in \cite{Speicher}, the set $\cup_{n>1} NC(n) \times \M^{\otm B}$ and a family of functions $\hat{\phi}=(\phi^{(n)})_{n>1}: \cup_{n>1} NC(n) \times \M^{\otm B} \rightarrow B$ and defined $\hat{\phi}(\pi)(a_1 \otm \ldots \otm a_n)$, where $\pi$ is a noncrossing partition, as explained in the previous section. The idea is to replace the object  $\cup_{n>1} NC(n) \times \M^{\otm B}$
by $\DM$ equipped with the $\NCP(\B)$-operad $\P[\hat{\phi}]$ and to reformulate a result of R. Speicher \cite{Speicher}. 
A \textit{moment function} $\phi$ is defined by $\phi^{(1)}=1$ and by $(n>1)$, 
$$\phi^{(n-1)}(a_1 \otm \ldots \otm a_pa_{p+1} \otm \ldots \otm a_n)=
\phi^{(n)}(a_1 \otm \ldots \otm a_p \otm a_{p+1} \otm \ldots \otm a_n).$$
In this case, one can choose $\phi^{(n)}(a_1 \otm \ldots \otm a_n):=\phi^{(n)}(a_1a_2 \ldots a_n)$.
In our
framework, 
R. Speicher showed also that the cumulant function $\hat{C}$ obtain by convolution of $\hat{\phi}$ with the $Zeta$ function associated with the reffinement order of the noncrossing partitions is still a map from $\DM$ to $B$.
We now reformulate the result of Speicher \cite{Speicher}.
\begin{theo}
Fix $B$, a unital associative algebra. Let $(\M, \phi)$, be a noncommutative probability space and $\phi$ a moment function. Let $\M_1, \ldots, \M_n,$ be $n$ unital $B-B$-subalgebras of $\M$. Fix  $\M_1, \ldots, \M_n,$ $B-B$-subalgebras of $\M$. Consider the set $I:=\{X^{\vec{n}}\otimes a_1 \otm \ldots \otm a_n \in \M; \forall n>1; \textrm{such that} \ \exists \ i,j \ a_i \in \M_{\epsilon_i}, \ a_j \in \M_{\epsilon_j} \ \textrm{and} \ \epsilon_i \not= \epsilon_j \}$. Then,
$\M_1, \ldots, \M_n,$ are stochastically free if and only if $I \subseteq \ker \hat{C}$, where $\hat{C}:\DM \rightarrow B$, is the cumulant function associated with $\hat{\phi}$ \textit{via} the convolution with the Zeta function with respect to the reffinement order.
\end{theo}
\section{Conclusion and open questions}
One of the main results of this paper is a reformulation of the free dendriform dialgebra on one generator in a tractable way, that is \textit{via} a natural coding of trees in terms of parentheses. With the identification of rooted planar binary trees with noncrossing partitions \textit{via} a `projection' method, we have implemented an action of trees
over tensors of a $B-B$-bimodule $\M$ and pointed out a connection with free probability in terms of the free dendriform dialgebra generated by $\treeA$. From this point of view, the suitable combinatorial object to deal with free probability would be planar rooted binary trees equipped with the Tamari and the reffinement partial orders. 

\noindent
What has been done at the level of rooted planar binary trees can be considered for planar rooted trees, a super Catalan object. The object $\DM$ becomes $\TrM$, the free dendriform trialgebra over $\M$ \cite{LodayRonco}, see also \cite{Lerlat}. Therefore, does there exist a `super Catalan analogue' of noncrossing partitions and  a `super Catalan analogue' of free probability?

\noindent
\textbf{Acknowledgments:}
The author would like to thank Michael Sch\"urmann, Uwe Franz, Rolf Gohm and Stefanie Zeidler for their very warm hospitality during his stay at the Institut f\"ur Mathematik und Informatik, Greifswald, Germany, where this paper has been written.
\bibliographystyle{plain}
\bibliography{These}

\end{document}